\documentstyle[12pt]{article}

\newcommand{\R}{{\bf R}}

\newcommand{\Sp}{{\rm Sp}}

\newcommand{\spl}{{\bf sp}}

\newcommand{\Rr}{{\cal R}}
\newcommand{\Oo}{{\cal O}}

\newcommand{\Co}{{{\cal C}\!{\it onj}}}

\newcommand{\eps}{{\varepsilon}}
\newcommand{\de}{{\delta}}
\newcommand{\De}{{\Delta}}
\newcommand{\ga}{{\gamma}}

\newcommand{\la}{{\lambda}}
\newcommand{\si}{{\sigma}}

\newcommand{\Bb}{{\cal B}}
\newcommand{\Cc}{{\cal C}}
\newcommand{\Dd}{{\cal D}}

\newcommand{\Nn}{{\cal N}}

\newcommand{\Uu}{{\cal U}}

\newcommand{\Ss}{{\cal S}}

\newcommand{\proof}[1]{\noindent{\bf Proof#1:\  }}

\newtheorem{theorem}{Theorem}[section]

\newtheorem{corollary}[theorem]{Corollary}
\newtheorem{definition}[theorem]{Definition}

\newtheorem{remark}[theorem]{Remark}
\newtheorem{lemma}[theorem]{Lemma}

\newtheorem{prop}[theorem]{Proposition}

\newtheorem{makingDe}{Lemma \ref{lemma: constructDe}}

\newtheorem{makingdifforder} {Lemma \ref{lemma: difforder}}

\input{psfig}

\catcode`\@=12

\def\Empty{}
\newcommand\oplabel[1]{
  \def\OpArg{#1} \ifx \OpArg\Empty {} \else
  	\label{#1}
  \fi}
		
\long\def\realfig#1#2#3#4{
\begin{figure}[htbp]
\centerline{\psfig{figure=#2,width=#4}}
\caption[#1]{#3}
\oplabel{#1}
\end{figure}}

\newcommand{\comm}[1]{}

\title{The Positive Fundamental Group of ${\Sp}(2)$ and ${\Sp}(4)$}
\date{November 12, 1997}
\author{Jennifer Slimowitz \\
         Depatment of Mathematics\\
         SUNY at Stony Brook\\
         Stony Brook, NY  11794\\
         jslimow@math.sunysb.edu}

\begin{document}

\maketitle

\abstract{In this paper, we examine the homotopy classes of positive loops in ${\Sp}(2)$ and ${\Sp}(4)$.  We show that two positive loops are homotopic
if and only if they are homotopic through positive loops.}

\section{Introduction}
A positive path in the group of real symplectic matrices ${\Sp}(2n)$ is a smooth path $A(t)$ whose derivative $A'$ satisfies
$$A'(t) = JP_t A(t)$$
where $P_t$ is a positive definite symmetric matrix (dependent on t) and $J$ is the standard complex structure. It is easy to see that positive paths are exactly those generated by negative definite time dependent quadratic Hamiltonians on ${\bf R}^{2n}$.  The simplest example of a positive path is the counter clockwise rotation $A(t) = e^{Jkt}$ where $k$ is any positive integer; here
$P_t = kI$. 

The relationship between positive paths and geodesics under the Hofer norm
motivates Lalonde and McDuff's discussion in ~\cite{LM2}.  In particular, a compactly supported Hamiltonian $H_t: {\bf R}^{2n} \rightarrow {\bf R}$ generates a flow $\psi _t$  for
$t \in [0,1]$ in the group of compactly supported Hamiltonian diffeomorphisms.
$\psi _t$ is a geodesic under the Hofer norm if and only if around each $t_0 \in [0,1]$ there exists
an interval $I$ such that there
exist two points $p$ and $P$ so that $p$ is a minimum and $P$ is a maximum of
$H_t$ for all $t \in I$.  Around $p$, the linearized flow of $\psi_t$ is
a positive path in ${\Sp}(2n)$, and it is called short if 1 is not an eigenvalue of the $\psi_t$  for any $t \in [0,1]$.   In \cite{BP} and
\cite{LM1}, it is shown that if the linearizations of the flow at $p$ and $P$ are short, then $\psi_t$ is a stable geodesic.   Lalonde and McDuff study this linearized 
flow and positive paths in general in order to obtain topological information about stable geodesics in \cite{LM2}.  Their work further develops 
Krein's theory by analyzing short positive paths whose 
whose eigenvalues lie off of the unit circle.  They show that any short positive path may be extended to a short positive path whose endpoint is diagonalizable with eigienvalues on $S^1$, and also that the space of short positive paths which
end at such a matrix is path connected \cite{LM2}.  

In addition, they define
the positive fundamental group ${\pi}_{1,pos} ({\Sp}(2n))$ to be the semigroup generated by positive loops with base point at the identity, where two loops are considered equivalent if one can be deformed to the other via a family of
positive loops.  In \cite{LM2}, Lalonde and McDuff pose the natural question:  ``Is the obvious
map
$${\pi}_{1,pos} ({\Sp(2n)}) \rightarrow {\pi}_1 ({\Sp}(2n))$$
injective?''  This paper provides the answer in the affirmative for ${\Sp(2)}$
and ${\Sp}(4)$.   The main difficulty in the four dimensional case is to show that any positive loop is  positively homotopic to a loop whose eigenvalues
lie in $S^1 \cup {\bf R}$.  

To examine the behavior of positive loops in ${\Sp}(2n)$, we follow the
lead of Lalonde and McDuff and look at the projection of these loops in 
the stratified space of symplectic conjugacy classes.  We characterize these
projections and construct homotopies between them,  and then lift the results to ${\Sp}(2n)$ by means of a lifting lemma.  Lalonde and McDuff look at generic
paths and those meeting isolated codimension two singularities;  here we will
occassionally need to look at paths which cross singularities of higher codimension. The notation in this paper is consistent with \cite{LM2}; their 
results will be quoted without proof.

I thank my advisor Dusa McDuff for introducing me to symplectic topology,
giving me insight into this problem, and commenting on many previous
drafts of this paper. 

\section{The Behavior of Positive Paths and Lifting Lemmas}
A useful tool for describing the movement of eigenvalues along a positive
path is the splitting number.  The notion of splitting number arises from Krein theory, described in \cite{EK1} and \cite {EK2}, and is explained further in Lalonde and McDuff ~\cite{LM2}.  They define the non-degenerate Hermitian symmetric form $\beta $ on ${\bf C}^{2n}$ by 
$\beta (v,w) = -i {\overline w} ^T J v$ where $J$ is the standard $2n \times 2n$
block matrix with the identity in the lower left box and minus the identity
in the upper right box. They prove the 

\begin{lemma}
If $A \in {\Sp}(2n)$ has eigenvector $v$ with eigenvalue $\lambda \in S^1$
of multiplicity 1, then $\beta (v,v) \in {\bf R} - \{0 \}$.
\end{lemma}

Hence, for any simple eigenvalue $\lambda \in S^1$ we may define $\sigma 
(\lambda) = \pm 1$ where $\beta (v,v) \in \sigma (\lambda) {\bf R} ^+$.
Using properties of $\beta$, we can check that $\sigma(\lambda) = - \sigma 
({\overline \lambda})$.  As an illustration, when $n = 1$, the  matrix
$$ J = \left  (\begin{array}{cc}  0 & -1 \\
                                  1 & 0 \end{array} \right)
$$
has eigenvalues $i$ and $-i$ corresponding to the eigenvectors $v_i =
\left (\begin{array}{c} 1 \\ -i  \end{array} \right)$ and $v_{-i} =
\left (\begin{array}{c} 1 \\ i  \end{array} \right)$.  Computing, we find
that $\beta (v_i, v_i) = 2$ so $\sigma(i) = 1$ and, similarly $\sigma(-i) = -1$.

In a more general setting, if $\lambda \in S^1$ has multiplicity $> 1$, we set
$\sigma(\lambda)$ to be equal to the signature of $\beta$ on the corresponding
eigenspace.  It is a straightforward calculation to see that the symplectic 
conjugacy class of a diagonalizable element in ${\Sp}(2n)$ with all of its eigenvalues on
the circle is determined by its spectrum and corresponding splitting numbers.
Hence, for each pair of conjugate eigenvalues  $ \{\lambda, {\overline \lambda} \} \in S^1$, there exist two symplectic conjugacy classes in ${\Sp}(2)$:  one where $\lambda$
has positive splitting number (and ${\overline \lambda}$ has negative
splitting number) and one where $\lambda$ has negative splitting number (and
${\overline \lambda}$ has positive splitting number).  Note that there is
no corresponding notion for real eigenvalues or the eigenvalues on $S^1$
of a non-diagonalizable matrix.  

A natural question to ask is, ``What restrictions does positivity impose 
upon movement of eigenvalues?"    Krein's lemma states that under a positive
flow, simple eigenvalues on $S^1$ with $+1$ splitting number move counter
clockwise while those with $-1$ must move clockwise \cite{EK2}.  In \cite{LM2}, Lalonde and McDuff
show that when a positive path has a pair of eigenvalues that enter $S^1$,
they must do so at a matrix which has a $2 \times 2$ Jordan block symplectically
conjugate to 
$$N_{\lambda}^+ =  \left(\begin{array}{cc}  {\lambda} &  -{\lambda} \\
                             0    &  {\lambda}    \end{array}\right)
$$
where ${\lambda}$ represents the eigenvalue on $S^1$.  Similarly, when a pair
leaves $S^1$, it does so at a matrix with a Jordan block symplectically
conjugate to 
$$N_{\lambda}^- =  \left(\begin{array}{cc}  {\lambda} &  {\lambda} \\
                             0    &  {\lambda}    \end{array}\right).
$$
These restrictions are, in fact, the only ones dictated on generic paths by the positivity condition, leaving us with the following statement:

\begin{lemma} \label{lemma: restrictions}
A positive path in ${\Co}$ may move freely between conjugacy classes when
its eigenvalues are away from $S^1$.  On $S^1$, the eigenvalues move according to splitting
number by
Krein's lemma, and when entering and leaving $S^1$, they behave according
to the above results of Lalonde and McDuff.
\end{lemma}

For example, there are $4$ open regions in ${\Sp}(4)$ whose union is dense: 

\begin{description}  
\item[(i)] ${\Oo}_{\Cc}$ , consisting of all matrices with 4 distinct eigenvalues of the form $\{ {\lambda}, \overline{\lambda},  \frac{1}{\lambda},
\frac{1}{\overline{\lambda}} \} \in \bf{C}-(\bf{R} \cup S^1)$; 
\item[(ii)] ${\Oo}_{\Uu}$ , consisting of all matrices with eigenvalues on $S^1 - \{1, -1 \}$
where each eigenvalue has multiplicity 1 or multiplicity 2 with non-zero splitting numbers;
\item [(iii)]  ${\Oo}_{\Rr}$ , consisting of all matrices whose eigenvalues have multiplicity 1 and lie on ${\bf R} - \{0,1,-1 \}$ 
\item[(iv)]  ${\Oo}_{\Uu,\Rr}$ , consisting of all matrices with 4 distinct eigenvalues, one pair on $S^1 - \{ 1, -1  \}$ and the other on ${\bf R} - \{0,1,-1\}$.
\end{description}
We will describe the other higher codimension regions later.  Lemma
\ref {lemma: restrictions} tells us that positive paths may move freely
in ${\Oo}_{\Cc}$ and ${\Oo}_{\Rr}$, but their behavior is restricted
when in  and when entering or leaving ${\Oo}_{\Uu}$ and ${\Oo}_{\Uu, \Rr}$.

Also useful will be the basic facts about positive paths from \cite{LM2}:

\begin {lemma} \label{lemma: openness}
\begin{description}
\item[(i)] The set of positive paths is open in the $C^1$ topology.
\item[(ii)] Any piecewise positive path may be $C^0$ approximated by a positive path. \end{description}
\end{lemma}

We now begin the discussion of homotopy and develop the tools necessary
to prove the injectivity of map from $\pi_{1, pos}({\Sp}(2n)) \rightarrow
\pi_1 ({\Sp}(2n))$.  Given a homotopy whose endpoints are positive paths,
we need to produce a homotopy between those two endpoints where each 
path in the homotopy is a positive path.  We will consider the projection
of the original homotopy to ${\Co}$, the space of symplectic conjugacy classes.
Let $\pi$ denote this projection:
$${\pi}(A) =  \cup_X \{ XAX^{-1} : X \in {\Sp}(2n) \} \in {\Co}. $$
After altering the projection of the homotopy in ${\Co}$ in a specific way to make each path in it
positive, we lift it to ${\Sp}(2)$ or ${\Sp}(4)$.  

Now we will state some useful definitions and two propositions which will
enable us to execute the lifting.

\begin{definition} 
A point in ${\Sp}(2n)$ is called a {\bf generic point} if all of its eigenvalues
have multiplicity 1.  
A path in ${\Sp}(2n)$ is called a {\bf generic path} if all of its points
are generic or lie on the codimension 1 boundary part of a generic region,
and the codimension 1 boundary points are isolated.    These definitions also
hold for points and paths in ${\Co}$.
\end{definition}

\begin{definition}
A path $a_t$ in ${\Co}$ is called {\bf positive} if there exists a positive
path ${A}_t \in {\Sp(2n)}$ such that $\pi ({A}_t) = a_t$.
A homotopy ${H}(s,t) \in {\Sp}(2n)$ is called {\bf positive} if for every $s_0$,
${H}(s_0,t)$ is a positive path.  A homotopy $h(s,t) \in {\Co}({\Sp}(2n))$ is called {\bf positive} if it is made up of positive paths in ${\Co}$, i.e. for every $s_0$, there is a positive path ${H}(s_0,t) \in {\Sp}(2n)$ such that
${\pi}({H}(s_0, t))  = h(s_0, t)$.
\end{definition}

\begin{prop}
Let ${A_t} \in {\Sp}(2n)$ be a generic positive path joining two generic points
${A_0}$ and ${A_1}$.   Then the set of positive paths  in ${\Sp}(2n)$ which 
lift $\pi ({A_t}) \in {\Co}$ is path connected.
\label{prop: Dusalifting}
\end{prop}

\proof{}
Here is the idea of Lalonde and McDuff's proof from \cite{LM2}.  Suppose ${B_t}$ and
${C_t}$ are two paths which lift $\pi( {A_t})$.  We may assume
that ${B_t}$ crosses codimension 1 strata at finitely many times $t_i$.
Note that each fiber of $\pi : {\Sp}(2n) \rightarrow {\Co}$ is path 
connected since ${\Sp}(2n)$ is.  
Hence, using Lemma \ref{lemma: openness}, we may homotop ${C_t}$ around those times  to $X_i {B_t} X_i^{-1}$ for $t$ close to $t_i$ for
some symplectic matrix $X_i$.  Let $\xi_B$ be the vector field tangent to
to the curve ${B_t}$ and define the vector field $\xi_C = X_i \xi_B X_i^{-1}$ over neighborhoods of each $ {C_{t_i}}$.  If we extend $\xi_C$  appropriately and take the convex combination vector fields $s \xi_B + (1-s) \xi_C$, these new
positive
vector fields have integral curves which also project to $\pi( {A_t})$.
Thus, the family of integral curves as $s$ varies from 0 to 1 gives a path
between  ${B_t}$ and ${C_t}$ within the lifts of $\pi( {A_t})$.
$\Box$

Certainly, if ${A_t}$ and ${B_t}$ are positively homotopic paths in ${\Sp}(2n)$,
then ${\pi}({A_t})$ and ${\pi}({B_t})$ are positively homotopic in ${\Co}$.  The
next proposition shows that when each path in the homotopy is generic, the converse is also true.

\begin{prop}
Let $h(s,t)$ be a positive homotopy of generic loops based at the identity in ${\Co}({\Sp}(2n))$ where
$ h:[0,1] \times [0, 2 {\pi} ] \rightarrow {\Co}$ and
$$ \begin{array} {rr}
h(0,t) = a_t &  h(1,t) = b_t  \end{array} . $$
Also, let ${A_t}, {B_t} :  [0, 2{\pi}] \rightarrow {\Sp}(2n)$ be any two positive generic loops based at $I$ so that  ${\pi}({A_t}) = a_t$ and 
${\pi}({B_t}) = b_t$.  Then, there exists a positive homotopy $
{H}(s,t) : [0,1] \times [0, 2 {\pi}] \rightarrow {\Sp}(2n) $ such that  ${H}(0,t) = {A_t}$ and ${H}(1,t)  = {B_t}$. \label{prop: homlifting}
\end{prop}

\proof{}
The proof of this proposition mimics that of the previous one, only here
we must introduce parameters.  After dealing with the technicalities of
locally lifting $h$ around each codimension 1 point as in Proposition
\ref{prop: Dusalifting}, we are left with a finite sequence of 
${H}^i(s,t): [s_i, s_{i+1}] \times [0, 2 {\pi}] \rightarrow {\Sp}(2n)$,
homotopies defined on some partition $[s_i, s_{i+1}]$ of $[0,1]$.  Here, 
${\pi}({H}^i(s,t)) = h(s,t)$ , and each loop ${H}^i(s_i,t) :
[0, 2 {\pi}] \rightarrow {\Sp}(2n)$ is a generic positive path.  Using 
Proposition \ref{prop: Dusalifting}, for each $i$, we glue ${H}^i(s_{i+1},t)$ to
${H}^{i+1}(s_{i+1},t)$ via a family of positive loops, all of which
project to $h(s_{i+1},t)$ in ${\Co}$, and let ${H}$ be the resultant
homotopy.  At the end of this paper, we give the full details concerning the
lifting of some specific homotopies in ${\Sp}(4)$.
$\Box$

Hence, to prove the injectivity of the map from $\pi_{1,pos} \rightarrow
\pi_1$, we need only construct a positive homotopy of generic paths in ${\Co}$ between
the projections of the two given endpoints.  This is exactly  what will happen
in ${\Sp}(2)$.   It turns
out, however, that the homotopy we construct  in ${\Co}({\Sp}(4))$
may have some paths which
are non-generic and go through points of codimension two and higher.
We will deal with this by finding specific lifts of the homotopy
in neighborhoods
of these points to ${\Sp}(4)$.   We then join these lifts
to the given homotopy using Proposition \ref {prop: Dusalifting}.

\section {The Positive Fundamental Group of ${\Sp}(2)$}

Here is a review of the structure of the stratified space of symplectic
conjugacy classes of ${\Sp}(2)$ as described in \cite{LM2}, along with some additional details. 

A generic matrix in ${\Sp}(2)$ has two distinct eigenvalues and belongs
to one of the following regions: 
\begin{description}
\item[(i)] ${\Oo}_{\Uu}$, consisting of all matrices with eigenvalues $\{        	\lambda, \overline{\lambda} \} \in S^1 $ 
\item[(ii)] ${\Oo}_{\Rr}$, consisting of all matrices with real eigenvalues $\{       \lambda, \frac{1}{\lambda} \}$ where $|{\lambda}| \geq 1 $.
\end{description}
We will divide each of ${\Oo}_{\Rr}$ and ${\Oo}_{\Uu}$ naturally into two parts: ${\Oo}_{\Rr}^+$ and
${\Oo}_{\Rr}^-$ for positive or negative eigenvalues and ${\Oo}_{\Uu}^+$ and
${\Oo}_{\Uu}^-$ based on the sign of the imaginary part of the eigenvalue with
positive splitting number. 
 
We see that the non-generic matrices are the identity matrix $I$ and
$-I$ and the non-diagonalizable matrices with a double eigenvalue of 1 or -1.
The space of symplectic conjugacy classes of ${\Sp}(2)$ (remember this requires similarity
by a symplectic matrix) can be described by the set $S^1 \cup (1, \infty) \cup
(- \infty, -1)$ in the plane with the points 1 and -1 tripled, as depicted in
Figure 1.  

 \realfig{fig1}{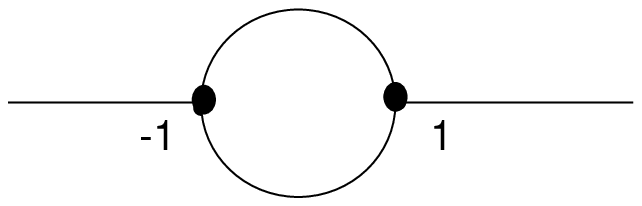}{${\Co}({\Sp}(2))$}{0.4\hsize}

This can be
seen as follows:  identify $A \in {\Oo}_{\Rr}$ with its eigenvalue ${\lambda}$
whose absolute value is greater than 1.  Clearly, all such matrices are conjugate.  For $A \in {\Oo}_{\Uu}$, we can distinguish between the two
eigenvalues $\{ \lambda, \overline{\lambda} \}$ by the notion of splitting
number as described above.  Associating $A$ to its eigenvalue with positive
splitting number produces a well-defined equivalence class, accounting for
each element in $S^1$.  $I$ and $-I$ each comprise their own equivalence class; associate $I$ with 1 and $-I$ with -1.  If $A$ is non-diagonalizable with double
eigenvalue -1, then A is conjugate to either 
$$
N_{-1}^+ = \left(\begin{array}{cc}  -1 &  1\\
                             0    &  -1    \end{array}\right)
$$
or 
$$
N_{-1}^- = \left(\begin{array}{cc}  -1 &  -1\\
                             0    &  -1    \end{array}\right).$$
In either case, we send $A$ to -1, and we have 3 conjugacy classes at
-1:  $-I, N_{-1}^+, N_{-1}^-$. 

Similarly, if $A$ is nondiagonalizable with double eigenvalue 1, then A is conjugate to either $N_{1}^+$ or $N_{1}^-$ , and we have
three conjugacy classes at 1:  $I, N_1^+, N_1^-$.  The space of $A \in {\Sp(2)}$
which project to either $N_{-1}^+,N_{-1}^-,N_1^+,$ or $N_1^-$ is of codimension 1.  By Lemma
\ref{lemma: restrictions} , we know that positive paths in ${\Co}$ enter $S^1$ 
via $N_{-1}^+$ and $N_1^+$ and leave via $N_{-1}^-$ and $N_1^-$.

\begin {definition}
A {\bf simple} path $\ga (t)$ in ${\Co}({\Sp}(2))$ has at most one local minimum and no local maxima each time it enters  $\pi (\Oo_{\Rr}^-)$ , and has at most one local maximum and no local minima each time it enters $\pi (\Oo_{\Rr}^+)$.

\end{definition}

\begin{lemma}
If ${\ga}_t$ is a simple path along the real axis in ${\Co}$ with bounded
eigenvalues, it is positive. \label{lemma: simple}
\end{lemma}

\proof{}
By Lemma \ref{lemma: openness}, it suffices to show that for all $\alpha, \beta \in
{\bf R}^+ - (0,1]$ where $\alpha \neq \beta$, there is a positive path ${A}_t \in
{\Sp}(2)$ such that the function $t \rightarrow \pi({A}_t)$ is an
embedding of $[0,1]$ onto $[\alpha,\beta]$ sending $0$ to $\alpha$ and
$1$ to $\beta$.  Consider the path $e^{Jt}B$ where 
$$B = \left ( \begin{array} {cc}
\beta & 0 \\
0 & \frac{1}{\beta}  \end{array}
\right ) $$
the projection of $e^{Jt}B$ to the real axis in ${\Co}$ depends only on the
trace of the matrix, since we can recover the eigenvalues from the trace
and the determinant which we know is 1.  So, by examining the movement of the
trace of $e^{Jt}B$, we can determine the flow of $\pi (e^{Jt}B)$ on the 
real axis.  We know that $\pi (e^{Jt}B)$ must travel counter clockwise
along the circle by Krein's Lemma, so once we figure out what the trace is
doing, we will get the trajectory of the path in all of ${\Co}$.
The derivative of the trace of $e^{Jt}B$ at time $t$ is 
$$ - (\beta + \frac{1}{\beta}) \sin t$$ 
which is negative for $0<t<\pi$, zero for $t=\pi$, and positive for $\pi<t<2 \pi$.   Note that at $t = \pi$, $e^{Jt}B = -B$.  Hence, $\pi(e^{Jt}B)$ 
finishes by coming off the circle through $N_1^-$ and travelling up the
real axis, past $\alpha$, to $\beta$.    We can let ${A}_t$ be the reparametrization of the last portion of $e^{Jt}B$ 
which projects to $[\alpha,\beta]$.   Similarly, if $\alpha > \beta$, 
we will let $A_t$ be the first part of $e^{Jt}C$ where 
$$C = \left ( \begin{array} {cc}
\alpha & 0 \\
0 & \frac{1}{\alpha} \end{array} \right ).$$
$\Box$

\begin {theorem}  
Suppose $A_t, B_t\in {\Sp}(2)$ are two positive loops based at $I$.  Then,
$A_t$ and $B_t$ are homotopic if and only if they are homotopic through
positive loops.  Thus, the natural map from
$$ \pi_{1,pos} ({\Sp}(2)) \rightarrow \pi_1 ({\Sp}(2)) $$
is injective and onto {\bf N}.  \label{theorem: mainthm2} \end{theorem}
\proof{}  
Certainly, if $A_t$ and $B_t$ are homotopic through positive loops, then
they are homotopic.

Conversely, if $A_t$ and $B_t$ are homotopic, then the homotopy descends to
a homotopy of the projections of the paths $\pi (A_t) $ and $\pi (B_t)$ in ${\Co}$.  Thus, the two projections
of the paths travel around $S^1$ the same number
of times; this homotopy invariant is the Maslov index.
 We can assume the paths only go though $I$ at times 0 and 1 and are generic
away from these points, as positivity is an open condition.  
We will show that $\pi (A_t)$ and $\pi (B_t)$ are both homotopic through positive
paths in ${\Co}$ to a standard path ${\ga _t}$ with appropriate Maslov index and, thus, that they
are homotopic through positive loops in ${\Co}$ to each other.  Since any piecewise positive path may be $C^0$ approximated by a positive path, it will
suffice to do the homotopy in pieces, first considering the parts of the
paths on the circle, and then considering the parts on the real axis.

Let ${\ga _t}$ be a loop at $I$ in ${\Co}$ which goes around $S^1$ the same 
number of times as $\pi (A_t)$ and $\pi (B_t)$.  Choose ${\ga _t}$ so that it is a simple path. (Thus, ${\ga _t}$
doesn't swivel back and forth more than once along the real axis each
time it leaves the circle.)  

Lemma  \ref{lemma: restrictions} tell us that by reparametrizing $\pi (A_t)$, we can make it equal to  ${\ga _t}$
for the times when ${\ga _t}$ takes values on the circle.  The new parametrization is positively homotopic to $\pi (A_t)$, so we need only
search for a homotopy from $\pi (A_t)$ to ${\ga_t}$ when these
paths take values on the real line.

From Lemma \ref{lemma: simple}, we know that the portion of each simple loop in ${\Co}$ on the real line is positive.  If $\pi (A_t)$ is simple, it can be easily homotoped through
positive paths to ${\ga _t}$ just by ``stretching'' through simple
and therefore positive paths.  If $\pi (A_t)$ is not simple, then we can
slightly perturb it to have finitely many local maxima and minima.  Then, we 
can consider each ``bump'' as a simple path, and flatten each one individually
by passing through simple and therefore positive paths.  After smoothing out all the bumps in this manner,  we are left
with one simple piecewise positive path in ${\Co}$ positively homotopic to
$\pi (A_t)$.  We can estimate this path closely by a simple postive path
positively homotopic to $\pi (A_t)$, and by 
moving through simple paths, homotop it to ${\ga _t}$.

Hence, $\pi (A_t)$ is homotopic through positive loops to ${\ga _t}$.
In the same way, $\pi (B_t) $ is also homotopic through positive loops to
${\ga _t}$, and so $\pi (A_t)$ and $\pi (B_t)$ are positively homotopic
in ${\Co}$.

All of these homotopies are through generic paths; hence by Proposition \ref{prop: homlifting}, we can lift this homotopy to ${\Sp}(2)$,
and the proof is complete. $\Box$

\begin{corollary}
Let $A_t : [0,2 \pi] \rightarrow  {\Sp}(2)$ be a positive loop.  Then $A_t$ is positively homotopic
to $e^{Jkt}K$ where $k$ is the Maslov index of $A_t$ and $A_0 = K$. 
\label{corollary: homtoeJkt}
\end{corollary}

\proof{}
Since the Maslov index completely dictates the homotopy class of a loop,
$A_t$ is homotopic to $e^{Jkt}K$.  Hence by Theorem \ref{theorem: mainthm2},
$A_t$ is positively homotopic to $e^{Jkt}K$. $\Box$

Here are a few interesting remarks concerning positive paths in ${\Sp}(2)$:
 
\begin{remark}
At a point $A \in \Sp(2)$,the intersection of the positive cone and the tangent vectors pointing in the direction of the conjugacy class of $A$ is $$
 \{ JPA \cap (MA-AM) \} $$ where $P$ is positive definite symmetric and $M \in \spl(2)$. If $$A = \left (\begin{array} {cc} \lambda & 0 \\
                                             0     & \frac{1}{\lambda}
   \end{array}\right) , {|\lambda|} > 1 $$
then this intersection is $$ \{ \left (\begin{array} {cc} 0& -z\lambda \\
                                        \frac{x}{\lambda} & 0  \end{array}\right) | \; x,z \in \R^+ \}. $$
The intersection of the positive cone and tangent vectors pointing within
the conjugacy class at $BAB^{-1}$ for $B \in \Sp(2) $ is 
$$
         \{ B\left (\begin{array} {cc}  0 & -z\lambda \\ 
                                \frac{x}{\lambda} & 0 \end{array}\right) B^{-1}
\; | \; x,z>0 \} .$$
Hence, if  $$M = \left ( \begin{array} {cc} a & b \\ c & -a \end{array} \right ),
 b,c > 0, $$
then the path $\ga(t) = e^{Mt}Ae^{-Mt}$ is  a 
positive path staying in the conjugacy class of $A$. 
\label{lemma: posconj}
\end{remark}

\begin{theorem}  \label{theorem: pathinclass}
Given any two elements in the same conjugacy class in $\Oo_{\Rr}^\pm \in \Sp(2)$, there exists a positive path within the conjugacy class from one to the other.
\end{theorem}
\proof{}
This is a direct result of Lobry's Theorem which may be stated as follows:  Let $M$ be a smooth, connected, paracompact manifold of dimension $n$ with a set of complete vector fields
$ \{X^i|i \in I\}$ for some index set $I$.  Consider the smallest family of vector fields containing the $X^i$ which
is closed under Jacobi bracket.  At each point of $M$, the values of the elements of this family are vectors in the tangent space to $M$ which generate a linear subspace ${S}$.    If $dim(S)=n$ for all points in $M$,
the positive orbit of a point under the vector fields $X^i$ is the 
whole manifold.  (See Lobry ~\cite{LOB}, Sussman ~\cite{SU}, and Grasse
and Sussman ~\cite{GS}.)

In our specific case, $M$ is the conjugacy class of an element
in $\Oo_{\Rr}^{\pm}$, a smooth, connected 2 dimensional paracompact manifold.
Let $A$ represent the diagonal element of this conjugacy class with eigenvalues $\lambda$ and $\frac{1}{\lambda}$.  Our index set
$I = \R^+ \times \R^+$ and our vector fields at $BAB^{-1}$ will be the positive vectors in
$T_{BAB^{-1}}M$ :
$$X^{x,z}_{BAB^{-1}} = B  \left (\begin{array} {cc} 0& -z\lambda \\
                                        \frac{x}{\lambda} & 0  \end{array}\right) B^{-1}.$$

At each point in $M$ the dimension of the subspace spanned by the $X^{x,z}$ is
2.  Closure under Jacobi bracket would only add more vector fields and hence increase the dimension of the subspace which is spanned, so the $dim(S) \geq 2$ at all points in the manifold.  But, $dim(S) \leq dim(T_{BAB^{-1}}M) = 2$, so
$dim(S) = 2$. Lobry's theorem applies, and we can move within the conjugacy
class positively from any one element to any other. $\Box$

\begin{remark}
There exist positive paths $\ga^{\pm}(t)\in \Sp(2) $ such that $$\lim _{t \rightarrow \infty} trace(\ga^{\pm}(t)) = \pm \infty.$$  
\end{remark}
\proof {}
It suffices to find $\ga ^+$, as then we could just set $-\ga ^+ = \ga^-$.
Take the path $\ga {^+} {_0} = e^{JPt}A_0$ where $$
A_0 = \left (\begin{array} {cc} \lambda_0 & 0 \\
                                             0     & \frac{1}{\lambda_0}
   \end{array}\right) , {\lambda_0} > 1 \; , 
P = \left (\begin{array} {cc} 1&-1 \\ 2&-1 \end{array} \right) $$
Note that $\ga {^+} {_0} (0) = A_0$.
If we take the derivative of the trace of $\ga {^+} {_0}$, we find that
$$
\frac{d}{dt} tr(e^{JPt}A_0) = \lambda_0(\cos t - \sin t) - \frac{1}{\lambda_0} (\sin t + \cos t) $$
which is positive for $t<\tan ^{-1} \left
(\frac {\lambda_0^2 - 1}{\lambda_0^2 +1}\right)$ and zero for
$t=\tan ^{-1} \left (\frac {\lambda_0^2 - 1}{\lambda_0^2 +1}\right)$. At this local maximum, the trace of  $\ga {^+} {_0}$ is $2\lambda_0^3 + \frac{2}{\lambda_0} > 2\lambda_0^3$.  

The idea for creating a positive path whose trace goes to $\infty$ involves
gluing together successive paths of the type $\ga {^+} {_0}$ using 
Lemma \ref{lemma: openness}.  We start at some
diagonal matrix $A_0$ as above and let the first leg of $\ga {^+}$ be $\ga {^+} {_0}$ until time $t_0=\tan ^{-1} \left (\frac {\lambda_0^2 - 1}{\lambda_0^2 +1}\right)$.  By Theorem \ref{theorem: pathinclass}, there exists a positive path in the conjugacy
class between $\ga {^+} {_0} (t_0)$ and the diagonal element representing this
conjugacy class, say $A_1$. We can glue this path and $\ga {^+} {_0}$ together to get a positive
path from $A_0$ ending at the diagonal element $A_1$ with $tr(A_1) > \lambda_0 ^3$.

We let the second leg of  $\ga {^+}$ be $\ga_1 {^+} (t) = e^{JPt}A_1$, or actually
some reparametrization of this path to obtain the part where trace increases
past $\lambda_0^9$ followed by a positive path in the conjugacy
class to the diagonal element $A_2$ with $tr(A_2)>\lambda_0 ^9$.  Continue
in this manner gluing paths together, using $e^{JPt}$ to increase the trace 
followed by a path to the diagonal element  of the conjugacy class. We can see that the resultant path
will have trace tending to $\infty$, as with each step the trace not only increases, but it grows in a polynomial fashion.     $\Box$

\section{Positive paths in ${\Sp}(4)$}
The next theorem is the four dimensional analog of Theorem \ref{theorem: mainthm2} .   ${\Co}(Sp(4))$ is substantially more
complicated than ${\Co}(Sp(2))$; here we briefly recall its topology as described in \cite{LM2}.  Remember, we have the splitting number 
we can associate to simple eigenvalues on the circle which gives us a notion
of directionality, but we have no corresponding
idea for other eigenvalues.

Generic regions:
\begin{description}  
\item[(i)] ${\Oo}_{\Cc}$ , consisting of all matrices with 4 distinct eigenvalues in  $\bf{C}-(\bf{R} \cup S^1)$; one conjugacy class for each quadruple;
\item[(ii)] ${\Oo}_{\Uu}$ , consisting of all matrices with eigenvalues on $S^1$
where each eigenvalue has multiplicity 1 (or multiplicity 2 with non-zero splitting numbers); four (or two) conjugacy classes for each quadruple corresponding to
the possible
splitting numbers;
\item [(iii)]  ${\Oo}_{\Rr}$ , consisting of all matrices whose eigenvalues have multiplicity 1 and lie on ${\bf R} - \{0,1,-1 \}$; one conjugacy class for each
quadruple;
\item[(iv)]  ${\Oo}_{\Uu,\Rr}$ , consisting of all matrices with 4 distinct eigenvalues, one pair on $S^1 - \{1,-1\}$ and the other on ${\bf R} - \{0,1,-1\}$; two conjugacy classes for each quadruple corresponding to the
possible splitting numbers of the pair on $S^1$.
\end{description}

Codimension 1 boundaries of these regions: 
\begin{description} 
\item[(v)] ${\Bb}_{\Uu}$ , consisting of all non-diagonalizable matrices whose spectrum
consists of a pair of conjugate points in $S^1 - \{1, -1\}$ each of multiplicity
2 and splitting number 0; two conjugacy classes for each quadruple:   ${\Bb}_{\Uu}^-$ containing those matrices from which positive paths enter 
${\Oo}_{\Cc}$ and ${\Bb}_{\Uu}^+$
containing those matrices from which positive paths enter ${\Oo}_{\Uu}$;
\item[(vi)] ${\Bb}_{\Rr}$, consisting of all non-diagonalizable matrices whose 
spectrum is a pair of distinct points $\lambda, 1/\lambda \in {\bf R} - \{0,1,-1
\}$ each of multiplicity 2; one conjugacy class for each quadruple;
\item[(vii)]  ${\Bb}_{\Uu, 1}$, consisting of all non-diagonalizable matrices
with spectrum $\{\lambda, \overline{\lambda} , \pm1,$ $\pm1 \}$ with
$\lambda \in S^1 - \{1, -1\}$; two conjugacy classes for each quadruple, corresponding
to  $N_1^-$ (call this one ${\Bb}_{\Uu, 1}^-$) and $N_1^+$ (call this one ${\Bb}_{\Uu, 1}^+$);
\item[(viii)] ${\Bb}_{\Rr, 1}$, consisting of all non-diagonalizable matrices
with spectrum $\{ \lambda, 1/\lambda ,$ $\pm1, \pm1 \}$ with $\lambda \in 
{\bf R} - \{0,1,-1\}$; two conjugacy classes for each quadruple, corresponding
to $N_1^-$ (call this one ${\Bb}_{\Rr, 1}^-$) and $N_1^+$ (call this one ${\Bb}_{\Rr, 1}^+$).
\end{description}

It is useful to remember that generic positive paths move from ${\Oo}_{\Uu}$ to ${\Oo}_{\Cc}$ through ${\Bb}_{\Uu}^-$ and move back into  ${\Oo}_{\Uu}$
through ${\Bb}_{\Uu}^+$.
Postive paths 
from  ${\Oo}_{\Rr}$ to ${\Oo}_{\Cc}$  and  from 
 ${\Oo}_{\Cc}$ to ${\Oo}_{\Rr}$ pass through ${\Bb}_{\Rr}$.  Positive paths  going
from ${\Oo}_{\Uu}$ to ${\Oo}_{\Uu, \Rr}$ pass through  ${\Bb}_{\Uu, 1}^-$ and they return to ${\Oo}_{\Uu}$ through ${\Bb}_{\Uu, 1}^+$.  Finally, positive paths 
moving from ${\Oo}_{\Uu, \Rr}$ to ${\Oo}_{\Rr}$ pass through  ${\Bb}_{\Rr, 1}^-$, and those returning to  ${\Oo}_{\Uu, \Rr}$ pass through
${\Bb}_{\Rr, 1}^+$.

In addition, there are two important strata of higher codimension:
\begin{description}
\item[(ix)] ${\Bb}_{\Rr , \Dd}$, consisting of all diagonalizable matrices
with two real eigenvalues each of multiplicity two; 1 conjugacy class for
each quadruple;
\item[(x)] ${\Bb}_{\Uu , \Dd}$, consisting of all diagonalizable matrices with
a conjugate pair of eigenvalues on $S^1$, each of multiplicity two with 0
splitting number;  1 conjugacy class for each quadruple.
\end{description}

We now state the main theorem of this paper:

\begin{theorem}  \label{theorem: mainthm4}
Let ${A_t} , {B_t} : [0, 2{\pi} ] \rightarrow {\Sp}(4)$ be positive loops in $\Sp (4)$ with base point $I$.   Then ${A_t}$ and ${B_t}$ are homotopic if and only if they are 
homotopic through positive loops.  Hence, the natural map 
$${\pi}_{1, pos}({\Sp}(4)) \rightarrow {\pi}_1({\Sp}(4))$$
is injective and onto ${\bf N} - \{1\} $.
\end{theorem}

Certainly, if two loops are homotopic through positive loops, then they are homotopic.

The proof of the converse will come in several steps.  By Proposition 
\ref{prop: homlifting}, it will be sufficient to produce the positive
homotopy of generic loops 
in ${\Co}$ which can be lifted to $\Sp (4)$.  We will carefully examine the stratification of ${\Co}$ to determine the behavior of a generic positive path.
The idea is to first show that ${\pi}({A_t})$ and ${\pi}({B_t})$ can be positively
homotoped out of ${\pi}({\Oo}_{\Cc})$, leaving two loops in ${\Co}$ postively homotopic to ${\pi}({A_t})$ and ${\pi}({B_t})$  which are entirely contained in 
$${\Ss} = {\pi}({\Oo}_{\Uu}) \cup {\pi}({\Oo}_{\Rr}) \cup {\pi}({\Oo}_{\Uu , \Rr}) \cup {\pi}({\Bb}_{\Uu, 1}) \cup {\pi}({\Bb}_{\Rr, 1}) \cup {\pi}({\Bb}_{\Uu, \Dd}).$$
${\Ss}$ is the set of all open strata with eigenvalues in $S^1 \cup {\bf R}$
along with some boundary components to make it a connected set.  Then, we view these paths as residing in ${\Co}(\Sp (2)) \times {\Co} (\Sp (2)) \subset {\Co}(\Sp (4))$, allowing us to use results about $\Sp (2)$.  Finally, we show that two standard paths
which have eigenvalues traversing the circle with different speeds but with
the same number of total rotations are positively homotopic.  Using these lemmas we 
produce the homotopy in ${\Co}$, and then lift it to ${\Sp}(4)$ to prove the theorem.  We will postpone the technical proofs to the last section.

\begin{lemma}
Let ${A_t}$ be a positive generic loop with base point $I$.  Then,
${\pi}({A_t})$ can be positively homotoped out of ${\pi}({\Oo}_{\Cc})$ to a loop 
contained in ${\Ss}$.  
\label{lemma: outofOc}
\end{lemma}
\proof{}

We can slightly perturb any path so that it enters ${\Oo}_{\Cc}$ only a finite
number of times, hence we assume that ${\pi}({A_t})$ enters ${\pi}({\Oo}_{\Cc})$ only a finite number of times.  Krein shows that the very beginning and end of positive loops based at the identity must be in ${\Oo}_{\Uu}$.  More specifically, he shows that there exist positive ${\epsilon}$ and ${\delta}$
such that for all times $t$ where $0<t<{\epsilon}$
and $2{\pi} - {\delta} < t < 2{\pi}$ the path is in ${\Oo}_{\Uu}$ \cite{EK2}.  Therefore we need to consider the different ways in which ${\pi}({A_t})$ can leave ${\pi}({\Oo}_{\Uu})$, enter
${\pi}({\Oo}_{\Cc})$, and return to ${\pi}({\Oo}_{\Uu})$, and construct positive homotopies from each type to paths in ${\Co}$ which remain in ${\Ss}$.  Then, we can positively homotop each escape into ${\pi}({\Oo}_{\Cc})$ back into ${\Ss}$ individually to result in a loop in ${\Co}$ postively homotopic to 
 ${\pi}({A_t})$ and entirely
contained in ${\Ss}$.

First, notice that no positive path can travel directly from ${\Oo}_{\Uu, \Rr}$ to ${\Oo}_{\Cc}$ or ${\Oo}_{\Cc}$ to  ${\Oo}_{\Uu, \Rr}$ without crossing a boundary component of codimension greater than one.  Therefore, since ${A_t}$ is generic, it cannot contain these transitions.  Similarly, 
${A_t}$ cannot go directly from ${\Oo}_{\Uu}$ to ${\Oo}_{\Rr}$ or vice versa without crossing a higher codimensional boundary;  to avoid this, it must pass through ${\Oo}_{\Uu, \Rr}$  or ${\Oo}_{\Cc}$ at an intermediate time.

If ${\pi}({A_t})$ travels directly from ${\pi}({\Oo}_{\Cc})$ to  ${\pi}({\Oo}_{\Rr})$ and back to ${\pi}({\Oo}_{\Cc})$, Lalonde and McDuff
show how it can be perturbed to stay only in  ${\pi}({\Oo}_{\Cc})$ ~\cite{LM2}. 

Taking this into account, there are four distinct ways for ${\pi}({A_t})$ to leave ${\pi}({\Oo}_{\Uu})$, enter
${\pi}({\Oo}_{\Cc})$, and return to ${\pi}({\Oo}_{\Uu})$:   
$$  \begin{array}{ll}
 {\pi}({\Oo}_{\Uu}) \Rightarrow {\pi}({\Oo}_{\Cc}) \Rightarrow  {\pi}({\Oo}_{\Uu}) & (1)\\
 {\pi}({\Oo}_{\Uu}) \Rightarrow {\pi}({\Oo}_{\Cc}) \Rightarrow  {\pi}({\Oo}_{\Rr}) \Rightarrow {\pi}({\Oo}_{\Uu, \Rr}) \Rightarrow  {\pi}({\Oo}_{\Uu}) & (2) \\
{\pi}({\Oo}_{\Uu}) \Rightarrow  {\pi}({\Oo}_{\Uu, \Rr}) \Rightarrow {\pi}({\Oo}_{\Rr}) \Rightarrow {\pi}({\Oo}_{\Cc}) \Rightarrow  {\pi}({\Oo}_{\Uu}) & (3) \\
{\pi}({\Oo}_{\Uu}) \Rightarrow  {\pi}({\Oo}_{\Uu, \Rr}) \Rightarrow {\pi}({\Oo}_{\Rr}) \Rightarrow {\pi}({\Oo}_{\Cc}) \Rightarrow {\pi}({\Oo}_{\Rr}) \Rightarrow {\pi}({\Oo}_{\Uu, \Rr}) \Rightarrow  {\pi}({\Oo}_{\Uu}) & (4)
\end {array}  $$

At each transition, the path crosses the appropriate codimension one boundary.  Note that in each case, when the path is in ${\pi}({\Oo}_{\Rr})$, it has either come directly from or will go directly into ${\pi}({\Oo}_{\Cc})$. 
When passing between ${\pi}({\Oo}_{\Rr})$ and ${\pi}({\Oo}_{\Cc})$, both real
eigenvalues of multiplicty two are positive, or both are negative.  It is impossible to travel in real numbers from positive to negative without going through zero, and no symplectic matrix has 0 for an eigenvalue.  Therefore, all four eigenvalues will remain positive or all will remain negative for the entire time that ${\pi}({A_t})$ is in ${\pi}({\Oo}_{\Rr})$. 

Any generic positive path in ${\Co}$ can be broken up into finitely many
sections which lie in ${\Ss}$ connected by parts of type (1), (2),
(3), or (4).  Note that in between each escape into $\pi({\Oo}_{\Cc})$,
while the path is in ${\Ss}$, there is a time when one pair of eigenvalues
is $\{1,1\}$ and a time where one pair is $\{  -1,-1\}$.   This is due
to Lemma ~\ref{lemma: restrictions} and the fact that eigenvalues with
positive and negative splitting number must meet on $S^1$ in order for the
path to cross $\pi({\Bb}_{\Uu}^-)$ and enter $\pi({\Oo}_{\Cc})$.  Hence,
the different journeys into  $\pi({\Oo}_{\Cc})$ are separated by time and will not overlap at all.  If we could show how to positively homotop any path
of type (1), (2), (3), or  (4) back into ${\Ss}$, we could start with the 
first diversion that occurs (with respect to time) of $\pi(A_t)$ into
$\pi({\Oo}_{\Cc})$, homotop it back into ${\Ss}$, continue in the same
way one at a time with subsequent diversions, and eventually end up with
a path contained entirely in ${\Ss}$ and positively homotopic to $\pi(A_t)$.
Thus, the proof of Lemma \ref{lemma: outofOc} is now reduced to showing that
any path of type (1), (2), (3),  or (4) in ${\Co}$ is positively homotopic
to a positive path which lies in ${\Ss}$.

Note that (2) and (3) are opposites.  If we can perturb case (3) properly , then
we can also perturb case (2) in a similar manner.  Thus, we will only
work out the details for cases (1), (3), and (4).

\begin{lemma}  \label{lemma: type(1)}
Any path $a_t$ of type (1) in ${\Co}$  is positively  homotopic to a  positive
path which lies in ${\Ss}$.
\end{lemma}
\proof{}
Using Lemma \ref{lemma: restrictions} we can see that all paths of type
(1) with the same endpoints in $\pi({\Oo}_{\Uu})$ are homotopic.  It is therefore sufficient to consider a model path of type (1) and produce the homotopy for this case.
 We assume the eigenvalues of $a_t$  remain  on one pair of conjugate
rays in $\pi({\Oo}_{\Cc})$, and that $a_t$ simply goes out along these rays to a point where the norm of the
largest eigenvalue is $k$ and comes back.    Denote the elements of ${\pi}({\Bb}_{\Uu}^-)$ 
and $\pi ({\Bb}_{\Uu}^+)$  where $a_t$ enters and leaves
$\pi ({\Oo}_{\Uu})$  as  $\pi (X^-)$ and $\pi (X^+)$, respectively.
 
We will find a continuous family of positive paths in ${\Sp}(4)$
which leave ${\Oo}_{\Uu}$ at $X^-$, go into ${\Oo}_{\Cc}$ along
the appropriate ray, return along that ray, and re-enter ${\Oo}_{\Uu}$ at
$X^+$.  These paths should travel successively less far into ${\Oo}_{\Cc}$, with
their limit not going into ${\Oo}_{\Cc}$ at all, but staying on ${\Oo}_{\Uu}$
and passing through $X \in {\Bb}_{\Uu, \Dd}$.   Then, the projection of these paths to
${\Co}$ gives us the homotopy required by the lemma.   Note that if we find this continuous
family of positive paths for one $X$, we can do so for any other $Y \in 
{\Bb}_{\Uu , \Dd}$ by multiplying by $X^{-1}Y$.  Hence, without loss of generality,
we can assume that
$$ X =  \left(
\begin{array} {cccc} 
0 & 0 & 1 & 0 \\
0 & 0 & 0 & 1 \\
-1 & 0 & 0 & 0 \\
0 & -1 & 0 & 0    \end{array}  \right)$$
which has eigenvalues $\{i,i,-i,-i\}$.

Consider the path in ${\Sp}(4)$
$$ \gamma_k (r) = e^{Jr} \left(
\begin{array} {cccc}
0 & 0 & k & 0 \\
0 & 0 & 0 & 1/k \\
-k & 0 & 0 & 0 \\
0 & -1/k & 0 & 0 \\
\end{array}  \right)$$
as $r$ varies in a neighborhood of 0.  
The eigenvalues of $\gamma_k$ travel around the  circle, leaving at $X^-$
when $r$ is such that 
 $$\cos ^2 r = \frac{4k^2}{(1+k^2)^2}$$
to travel up the imaginary axis to the point $\{ ki, -ki, i/k, -i/k \} =
\pi (\gamma_k (0))$.  
Then, they move back down the imaginary axis to $X^+$, and re-enter the circle.
All the while that $\gamma_k$ is in ${\Oo}_{\Cc}$, its eigenvalues stay
on the imaginary axis.  The family $\gamma_k$ as we let $k \rightarrow 1$
is the continuous family of positive paths we need.  Note that the last path
in the homotopy will go through the non-generic stratum ${\Bb}_{\Uu, \Dd}$.
$\Box$

Case (4) requires us to consider exactly what part of ${\pi}({\Oo}_{\Rr})$ ${\pi}({A_t})$ enters.  First consider the case where both journeys into ${\pi}({\Oo}_{\Rr})$  are in ${\pi}({\Oo}_{\Rr}^+)$ or both are in ${\pi}({\Oo}_{\Rr}^-)$.  We know from Lemma \ref{lemma: restrictions} that movement in
${\pi}({\Oo}_{\Cc})$ is unrestricted by the positivity condition;
hence we can positively collapse the portion in ${\pi}({\Oo}_{\Cc})$ back to
either ${\pi}({\Oo}_{\Rr}^+)$ or ${\pi}({\Oo}_{\Rr}^-)$.  If, instead, this 
part of ${\pi}(A_t)$ moves ${\pi}({\Oo}_{\Uu}) \Rightarrow  {\pi}({\Oo}_{\Uu, \Rr}) \Rightarrow {\pi}({\Oo}_{\Rr}^+) \Rightarrow {\pi}({\Oo}_{\Cc}) \Rightarrow {\pi}({\Oo}_{\Rr}^-) \Rightarrow {\pi}({\Oo}_{\Uu, \Rr}) \Rightarrow  {\pi}({\Oo}_{\Uu})$ or  its opposite,  the analysis is more complicated.  We will call these cases (4a) and come back to them later.

Now let us consider case (3).   Assume without loss of generality that ${\pi}({A_t})$ enters ${\pi}({\Oo}_{\Rr}^+)$ instead of ${\pi}({\Oo}_{\Rr}^-)$.  
We can describe scenario (3) by graphing the motion of the eigenvalues in the complex plane as in Figure 2.

\begin{figure}[htbp]

\centerline{\psfig{figure=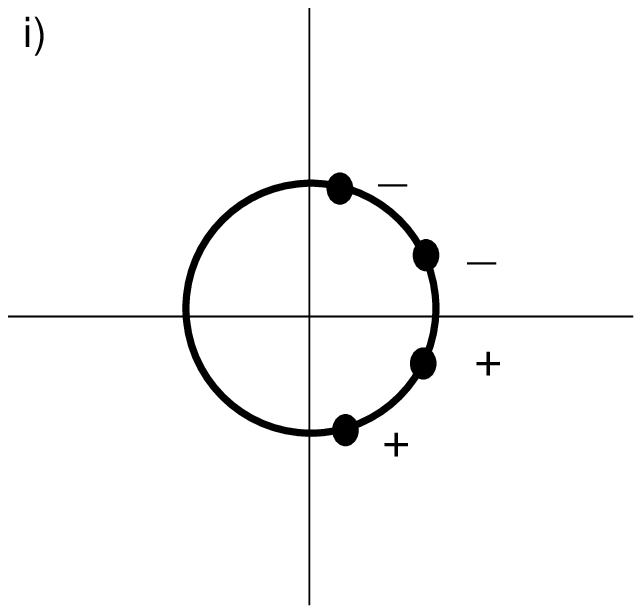,width=0.4\hsize}
\psfig{figure=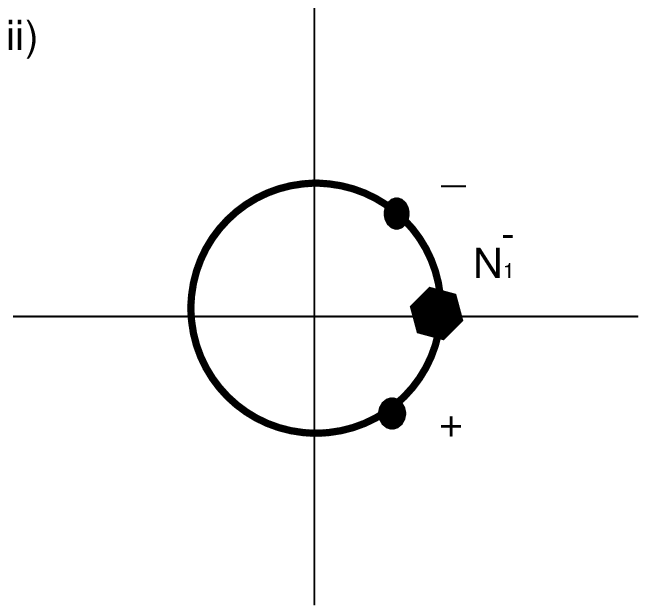,width=0.4\hsize}
\psfig{figure=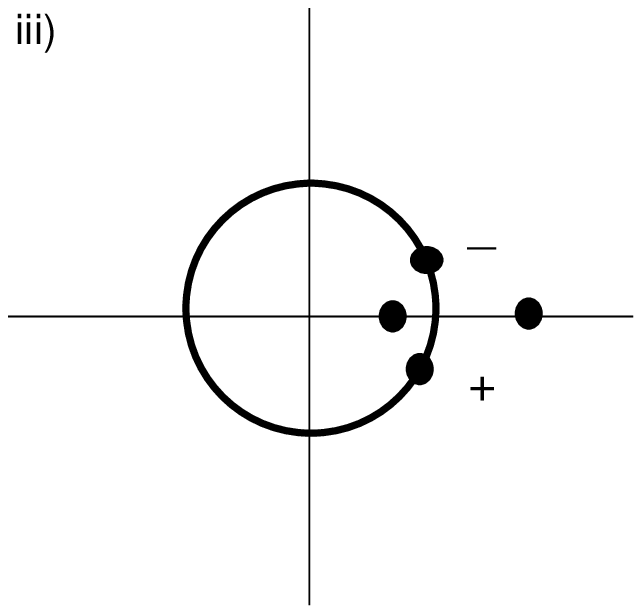,width=0.4\hsize}}

\vspace{1.5cm}

\centerline{\psfig{figure=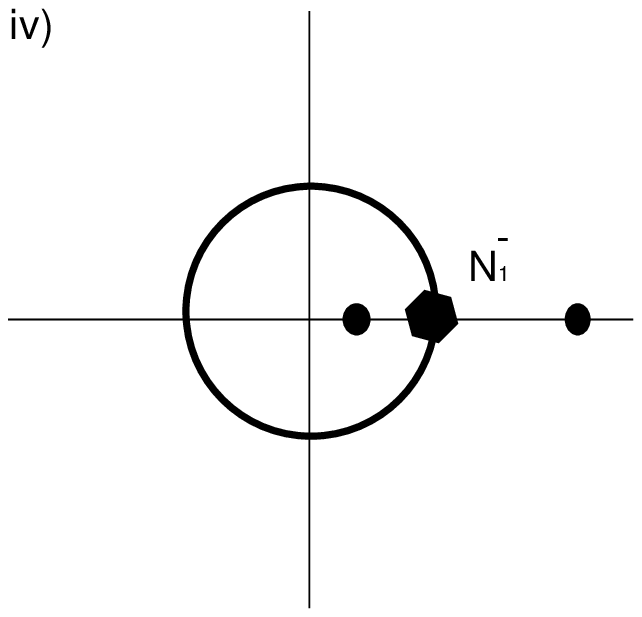,width=0.4\hsize}
\psfig{figure=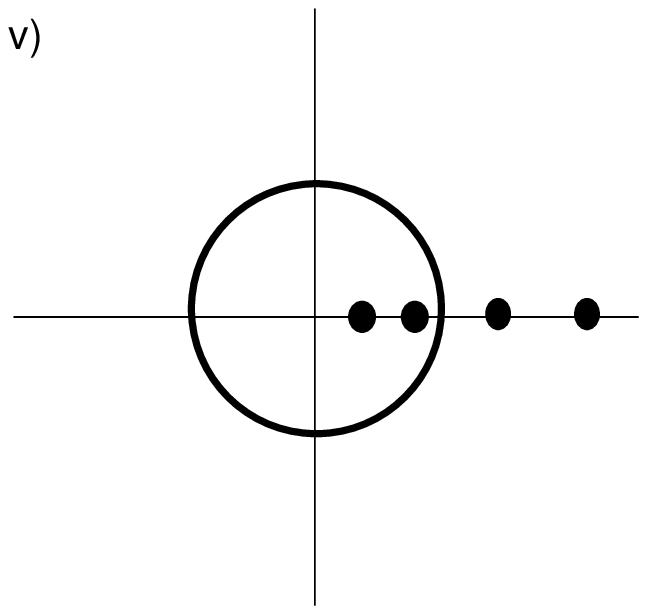,width=0.4\hsize}
\psfig{figure=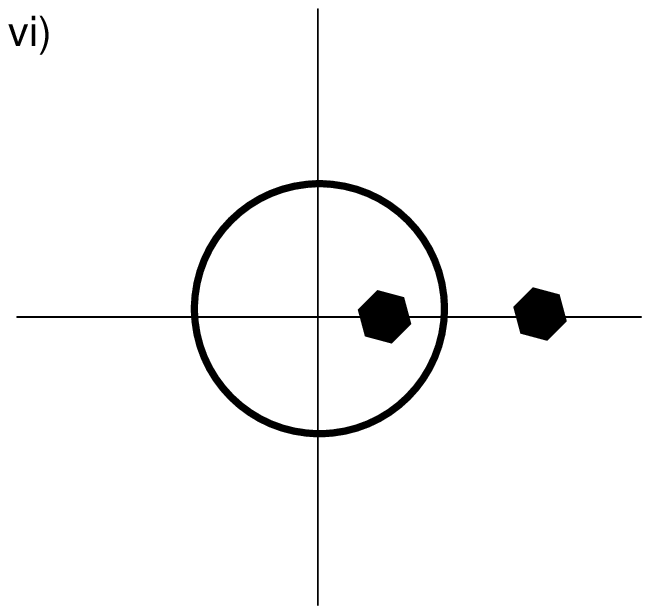,width=0.4\hsize}}

\vspace{1.5cm}

\centerline{\psfig{figure=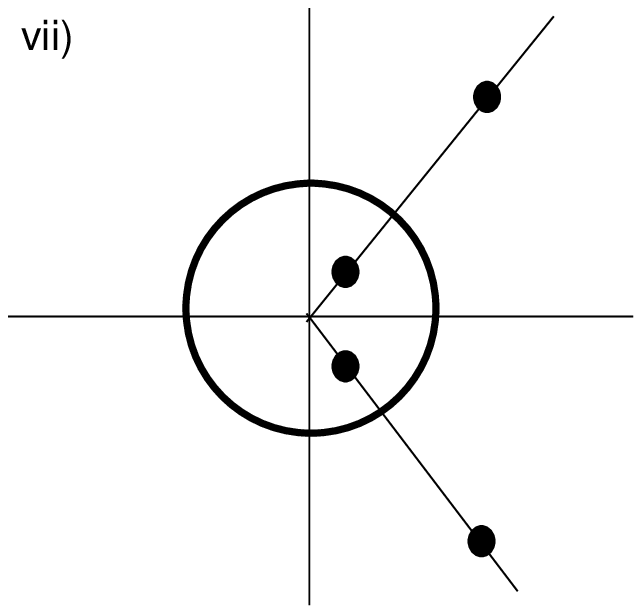,width=0.4\hsize}
\psfig{figure=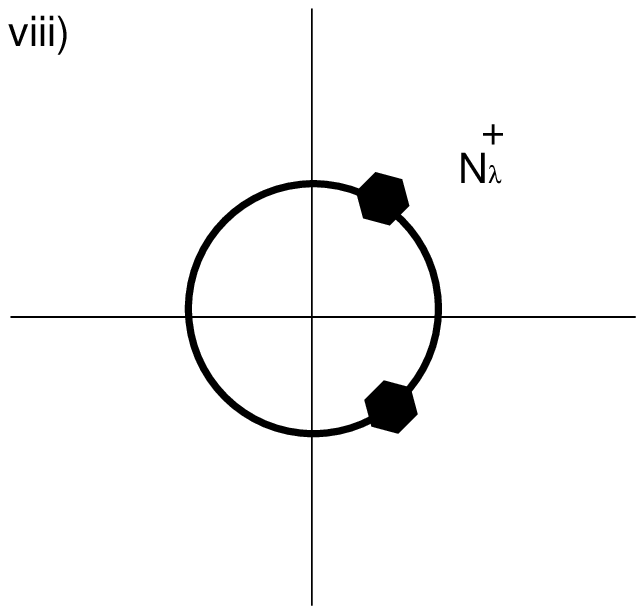,width=0.4\hsize}
\psfig{figure=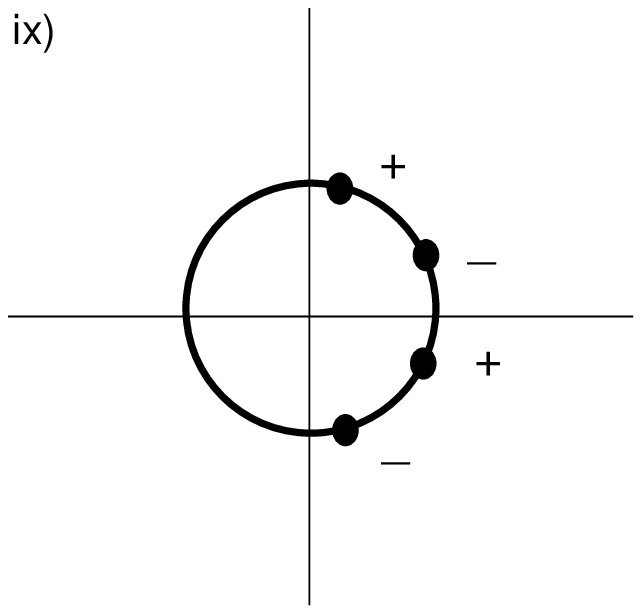,width=0.4\hsize}}

\caption[fig2]{$a_t$ in Case 3}
\oplabel{fig2}
\end{figure}

To begin, all four eigenvalues are on the circle, two conjugate pairs approaching the real axis.  Then, the first pair passes through ${\Nn}_1^-$, enters the real axis and the path is in ${\pi}({\Oo}_{\Uu, \Rr})$.  The second pair, still on the circle, migrates to the real axis also, eventually meets the first pair, and we have two real eigenvalues of multiplicity two.  At that moment, which we assume to be $t = \frac{1}{2}$, the path
breaks into ${\pi}({\Oo}_{\Cc})$.  Eventually the eigenvalues return to the 
circle as two plus/minus pairs, and continue rotating in the required direction.

\begin{lemma}  \label{lemma: type(3)}
Any path $a_t$ of type (3) in ${\Co}$ is positively homotopic to a positive
path which is contained in ${\Ss}$.
\end{lemma}
\proof{}
Using Lemma \ref{lemma: restrictions} it is easy to see that all generic paths
of type (3) with the same end points in $\pi({\Oo}_{\Uu})$ are positively homotopic.  Therefore, it suffices to start with
one path of this type and first show how to homotop it to a certain standard path $b_t$.  $b_t$ 
has the same first two configurations as  $a_t$, but, instead of the second pair entering
the real axis, the first pair re-enters the circle and the path is in ${\pi}({\Oo}_{\Uu})$ again.   Then, the positive eigenvalue from the first pair meets the eigenvalue with negative splitting number from the second pair, and vice versa, and the path escapes into ${\pi}({\Oo}_{\Cc})$.  Finally, this path returns directly to ${\pi}({\Oo}_{\Uu})$.   $b_t$ is depicted in Figure 3.

\begin{figure}[htbp]

\centerline{\psfig{figure=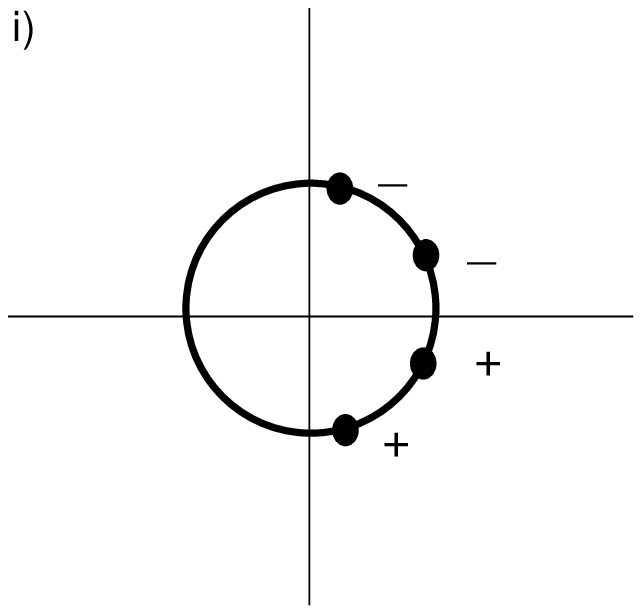,width=0.4\hsize}
\psfig{figure=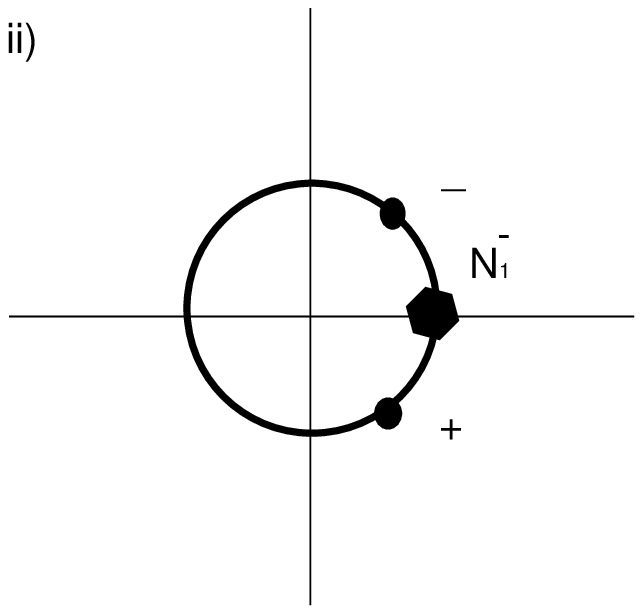,width=0.4\hsize}
\psfig{figure=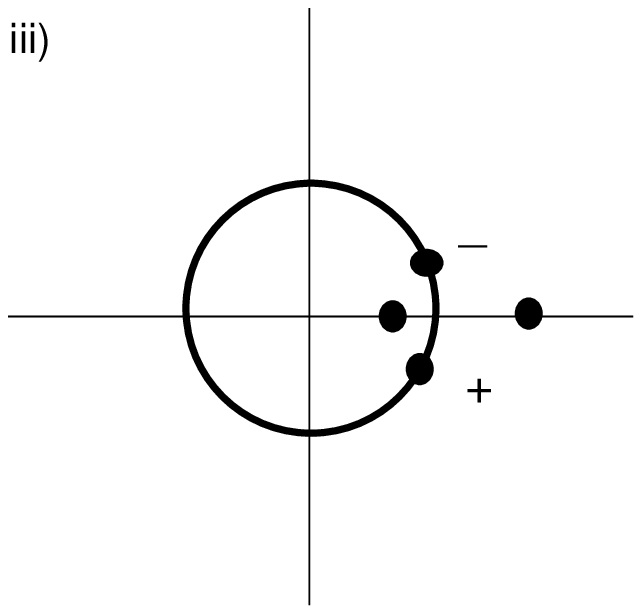,width=0.4\hsize}}

\vspace{1.5cm}

\centerline{\psfig{figure=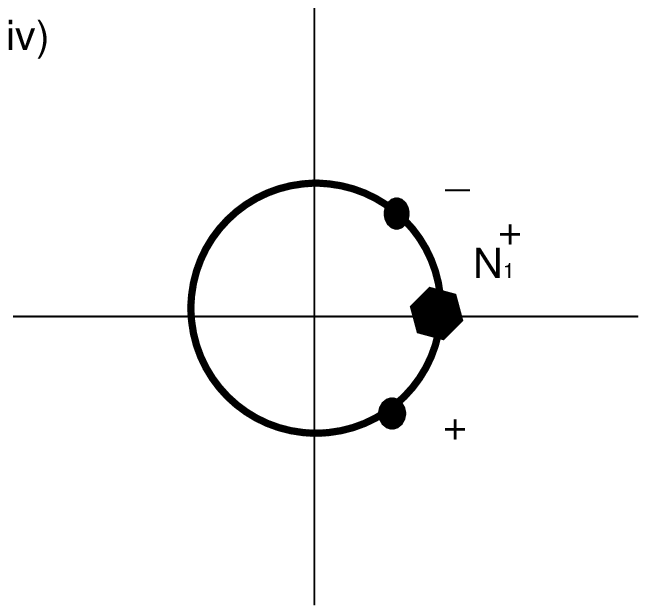,width=0.4\hsize}
\psfig{figure=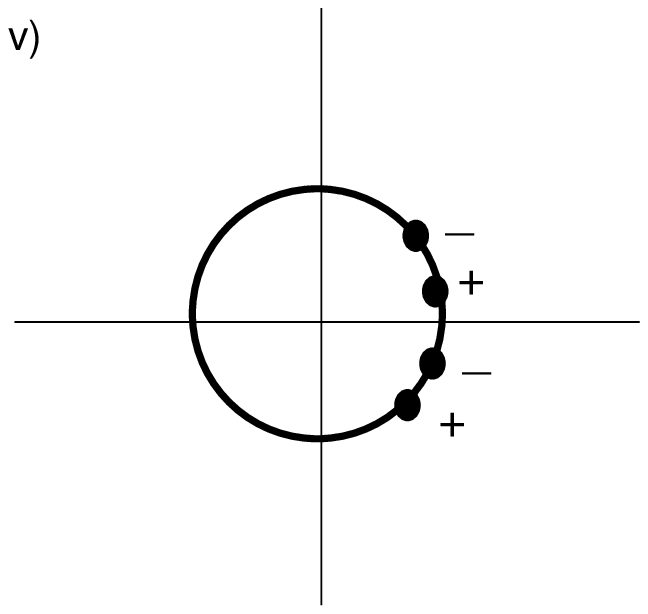,width=0.4\hsize}
\psfig{figure=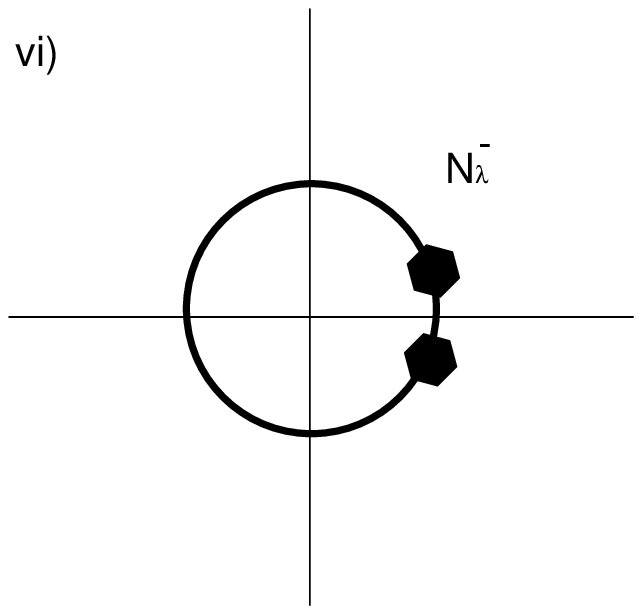,width=0.4\hsize}}

\vspace{1.5cm}

\centerline{\psfig{figure=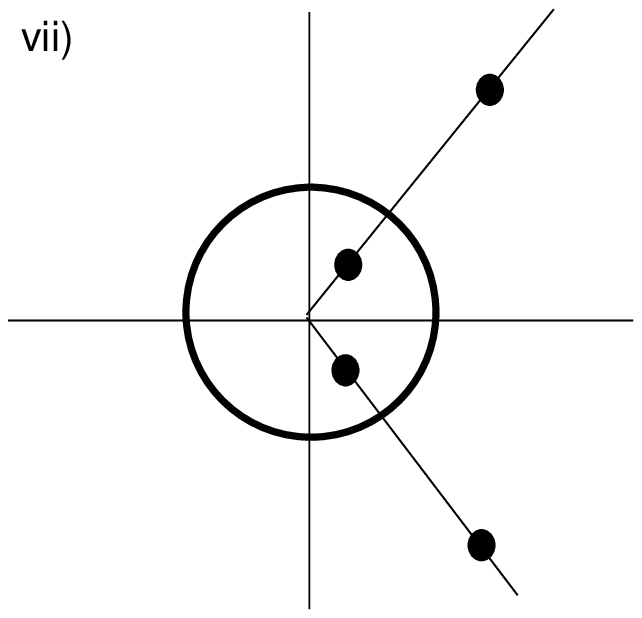,width=0.4\hsize}
\psfig{figure=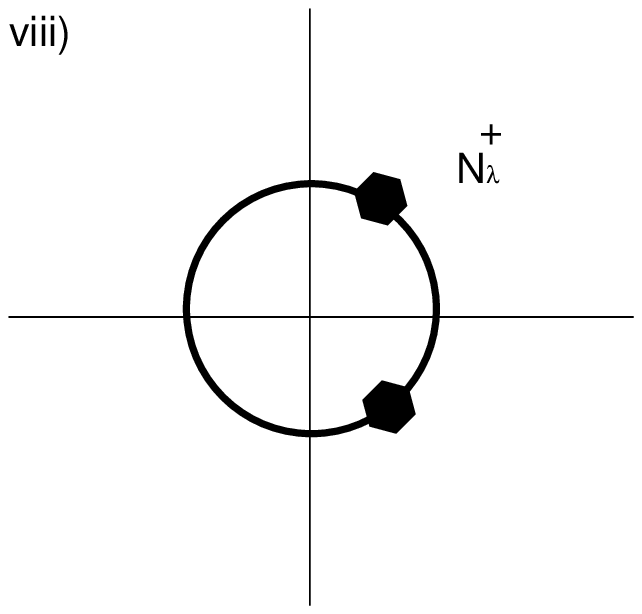,width=0.4\hsize}
\psfig{figure=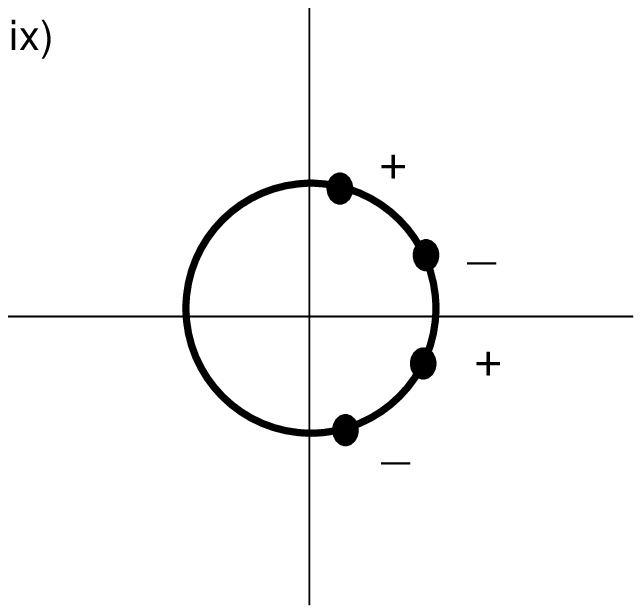,width=0.4\hsize}}

\caption[fig3]{The standard path $b_t$}
\oplabel{fig3}
\end{figure}

The path $b_t$ is of type (1), and we have already shown in Lemma
\ref{lemma: type(1)} how to positively homotope such paths out of ${\pi}({\Oo}_{\Cc})$.  Therefore, if we can construct the positive homotopy from  $a_t$ to ${b_t}$, we will be done with case (3).

The family of positive paths $a_t^s$ (where $a_t^0 = a_t$ and $a_t^1 = b_t$) we need to construct will all start in ${\pi}({\Oo}_{\Uu})$ and then go into ${\pi}({\Oo}_{\Uu, \Rr})$.  Here, $s$ is the homotopy variable and $t$ is
the time variable.  The first paths in the family will then enter  ${\pi}({\Oo}_{\Rr}^+)$ and break away into
${\pi}({\Oo}_{\Cc})$ at time $t = \frac{1}{2}$, just as $a_t$ does.  The point at which the $a^s$
enter ${\pi}({\Oo}_{\Cc})$ will progressively get closer and closer to and eventually hit the class of some matrix with eigenvalues $\{1,1,1,1\}$ at $s = 
\frac{1}{2}$
in ${\Co}$.  The paths subsequent to this will not enter  ${\pi}({\Oo}_{\Rr}^+)$, but rather will go back to ${\pi}({\Oo}_{\Uu})$ from ${\pi}({\Oo}_{\Uu, \Rr})$.  These paths will enter ${\pi}({\Oo}_{\Cc})$ from  ${\pi}({\Bb}_{\Uu}^-)$ at 
time  $t = \frac{1}{2}$ at points starting from the class of the matrix
with 1 as a quadruple eigenvalue, and travel up the circle.  Every path in
the family will reach ${\pi}({\Oo}_{\Cc})$ and travel back
to ${\pi}({\Oo}_{\Uu})$ ending at the same point as  $a_t$ and $b_t$.

Since movement in  ${\pi}({\Oo}_{\Cc})$ is not restricted under positivity, it suffices to find the family of positive paths $a_s^t$ at the infinitesmal level.  We will only construct the path $a_{\frac{1}{2}} ^s$ and its forward and backward tangent vectors to $a_s^t$ at $t = \frac{1}{2}$, since the rest of the
construction is straightforward.  We need to find two continuous vector fields along a continuous (not necessarily positive) path 
$q^s = a_{\frac{1}{2}}^s \in \pi ({\Bb}_{\Rr} \cup {\Bb}_{\Uu}^-)$ (except
when $s = \frac{1}{2}$)  which goes from  $q^0 \in {\pi}({\Bb}_{\Rr}) $, through some point with eigenvalues $\{1,1,1,1\}$ at $s= \frac{1}{2}$, to  $q^1 \in {\pi}({\Bb}_{\Uu}^-)$. Here, $q^0$ is the point in ${\Co}$ where $a_t$ enters ${\pi}({\Oo}_{\Cc})$, and $q^1$ is the point in ${\Co}$
where $b_t$ enters ${\pi}({\Oo}_{\Cc})$.  We need to find one positive continuous vector field pointing into  ${\pi}({\Oo}_{\Cc})$ at every point along $q^s$, and one negative continuous vector field pointing into ${\pi}({\Oo}_{\Rr}^+)$ at the points on $q^s$ with real eigenvalues and pointing into  ${\pi}({\Oo}_{\Uu})$ for all other points on $q^s$. We will explicitly find a lift $Q^s$ of such a path and vector fields in ${\Sp}(4)$; their
projections to ${\Co}$ will satisfy the required properties.  We set
$${\Nn}_1^{-,-} =  \left( \begin{array}{cccc}
1 & 1 & 0 & 0 \\
0 & 1 & 0 & 0 \\
0 & 0 & 1 & 1 \\
0 & 0 & 0 & 1 \end{array} \right ). $$
The proof of Lemma \ref{lemma: type(3)} is now reduced to the following:

\begin{lemma}
\label{lemma: constructDe}
 There exists a path $Q ^s : [0,1] \rightarrow  {\Bb}_{\Rr} \cup {\Bb}_{\Uu}^- \cup {\Nn}_1^{-,-}$ where $Q^0 \in {\Bb}_{\Rr}$, $Q^{\frac {1}{2}}=
\Nn_1^{-,-}$, and $Q^1 \in {\Bb}_{\Uu}^-$ satisfying two properties:
\begin{description}
\item[(i)] There exists a (positive) vector field along $Q ^s$ pointing into
${\Oo}_{\Cc}$  of the form $JP{Q}^s$ for a positive definite $P$.
\item[(ii)] There exists a negative vector field along $Q ^s$ pointing
into ${\Oo}_{\Rr}^+$ when $Q ^s \in {\Bb}_{\Rr}$ and pointing into
${\Oo}_{\Uu}$ elsewhere.  
\end{description}

\end{lemma}
\proof{}
The proof of this lemma will be deferred to Section 5. $\Box$

Now, the analysis for case (3) is finished.  We leave case (2) to the reader because it is very similar to case (3), and we are left only with case (4a):   ${\pi}({\Oo}_{\Uu}) \Rightarrow  {\pi}({\Oo}_{\Uu, \Rr}) \Rightarrow {\pi}({\Oo}_{\Rr}^+) \Rightarrow {\pi}({\Oo}_{\Cc}) \Rightarrow {\pi}({\Oo}_{\Rr}^-) \Rightarrow {\pi}({\Oo}_{\Uu, \Rr}) \Rightarrow  {\pi}({\Oo}_{\Uu})$ or ${\pi}({\Oo}_{\Uu}) \Rightarrow  {\pi}({\Oo}_{\Uu, \Rr}) \Rightarrow {\pi}({\Oo}_{\Rr}^-) \Rightarrow {\pi}({\Oo}_{\Cc}) \Rightarrow {\pi}({\Oo}_{\Rr}^+) \Rightarrow {\pi}({\Oo}_{\Uu, \Rr}) \Rightarrow  {\pi}({\Oo}_{\Uu})$.  This case is a
combination of cases (2) and (3).  Homotop the first part of the path from
${\pi}({\Oo}_{\Uu}) \Rightarrow  {\pi}({\Oo}_{\Uu, \Rr}) \Rightarrow {\pi}({\Oo}_{\Rr}) \Rightarrow {\pi}({\Oo}_{\Cc})$ to  ${\pi}({\Oo}_{\Uu}) \Rightarrow  {\pi}({\Oo}_{\Uu, \Rr}) \Rightarrow {\pi}({\Oo}_{\Uu}) \Rightarrow {\pi}({\Oo}_{\Cc})$ exactly the same way as in case (3).  Then, homotop the
second part from ${\pi}({\Oo}_{\Cc}) \Rightarrow {\pi}({\Oo}_{\Rr}) \Rightarrow {\pi}({\Oo}_{\Uu, \Rr}) \Rightarrow  {\pi}({\Oo}_{\Uu})$ to ${\pi}({\Oo}_{\Cc}) \Rightarrow {\pi}({\Oo}_{\Uu}) \Rightarrow  {\pi}({\Oo}_{\Uu, \Rr}) \Rightarrow
{\pi}({\Oo}_{\Uu})$ exactly the same way as in case (2).  This leaves a type (1) path
in ${\Co}$ positively homotopic to ${\pi}(A_t)$ which travels  ${\pi}({\Oo}_{\Uu}) \Rightarrow  {\pi}({\Oo}_{\Uu, \Rr}) \Rightarrow {\pi}({\Oo}_{\Uu}) \Rightarrow {\pi}({\Oo}_{\Cc}) \Rightarrow {\pi}({\Oo}_{\Uu}) \Rightarrow  {\pi}({\Oo}_{\Uu, \Rr}) \Rightarrow {\pi}({\Oo}_{\Uu})$.  Since 
type (1) cases have already been examined, the proof of Lemma \ref{lemma: outofOc} is now complete. $\Box$

\begin{lemma}
If $a_t$ is a positive loop in ${\Ss}$ based at $I$ constructed  by the methods
of Lemma \ref{lemma: outofOc}, then $a_t$ is positively homotopic in ${\Co}$ to $$
{\pi} \left ( \begin{array} {rr}
e^{Jkt} & 0 \\
0 & e^{J{\ell}t}  \end{array} \right ) $$
for some positive integers $k,\ell$.  \label{lemma: 4x4std}
\end{lemma}

\proof{}
 Note that the only time $a_t$ may go through a point with two eigenvalues
of multiplicity two is when forced to go through $\Bb_{\Uu, \Dd}$ as in
Lemma \ref{lemma: type(1)}. However, for each double pair of eigenvalues,
there is only one conjugacy class in $\pi ({\Bb}_{\Uu,\Dd})$ which we can
write as the class of an element in block diagonal form.  Hence, we can find positive loops ${X_t}, {Y_t} \in {\Sp}(2) $ such that $$
{\pi} \left ( \begin{array} {rr}
{X_t} & 0 \\
0 & {Y_t} \end{array} \right ) = a_t \in {\Co}  $$
where $X_0 = Y_0 = I$.
The result now follows from Theorem \ref{theorem: mainthm2} $\Box $

The next lemma shows that the positive homotopy class of $$
 \left ( \begin{array} {rr}
e^{Jkt} & 0 \\
0 & e^{J{\ell}t}  \end{array} \right ) \in {\Sp}(4) $$
where $k, \ell > 0$ depends only on the sum $k+{\ell}$, the same invariant as the regular homotopy class.  Let $\sim_+$ mean positively homotopic.

\begin{lemma}
Let $k,\ell \geq 1$, $n \geq 3$, $n>k$, $n>\ell$ where $k,\ell,n$ are all integers.
Then, $$
\left ( \begin{array} {rr}
e^{Jkt} & 0 \\
0 & e^{J(n-k)t} \end{array} \right ) \sim_+
\left ( \begin{array} {rr}
e^{J\ell t} & 0 \\
0 & e^{J(n- \ell )t} \end{array} \right ) $$
where $J \in {\Sp}(2)$.  \label{lemma: difforder}
\end{lemma}

\proof{}
The straightforward but detailed proof of this lemma is deferred to Section 5.

Now we have all the tools to prove the main theorem.

\proof{ of Theorem \ref{theorem: mainthm4}}

Take ${A_t}$ and ${B_t}$ to be two
homotopic positive loops in ${\Sp}(4)$ based at $I$.   By Lemmas \ref{lemma: outofOc} , \ref{lemma: 4x4std} , and \ref{lemma: difforder} , ${\pi}({A_t})$
is positively homotopic to ${\pi}({B_t})$.   Denote this homotopy in ${\Co}$
by $h(s,t)$.

The final step in the proof will be to use $h(s,t)$ to produce a homotopy ${H}(s,t) \in {\Sp}(4)$ where ${H}(0,t) = a_t$ and ${H}(1,t) =
b_t$.  If all of the loops in $h(s,t)$ are generic in ${\Co}$ except at the base point $I$, then
by Proposition \ref{prop: homlifting}, $h(s,t)$ can be lifted to ${\Sp}(4)$ and the proof of the theorem is 
complete.  Consider the case, then, when some loop in $h(s,t)$ is not generic; i.e.
there exists some $s_0 \in [0,1]$ such that $h(s_0,t)$ passes through a boundary
component of codimension greater than 1 or stays in a codimension 1 boundary
stratum for more than one instant.  Note that ${\pi}({A_t})$ and ${\pi}({B_t})$  are generic, so one of
the steps in the construction of $h(s,t)$ above must have introduced this
nongeneric behavior.
There are three isolated ways in which this can happen:  

\begin{description}
\item[(i)]  by the construction in Lemma \ref{lemma: type(1)} where a path goes through the stratum of diagonalizable
elements with 2 pairs of double eigenvalues $\{ \lambda, \lambda, \overline{\lambda}, \overline{\lambda} \} $ on the circle, $\pi({\Bb}_{\Uu , \Dd})$,
\item[(ii)]  while being homotoped out of $\pi ({\Oo}_{\Cc})$, by the construction in case (2),(3)
or (4a) of Lemma \ref{lemma: outofOc} where the paths are forced to go through ${\pi}({\Bb}_{\Rr, \Dd})$ or $N_1^{-,-}$,
\item[(iii)]  in the proof of Lemma \ref{lemma: 4x4std}  where loops are forced to
pass through $I$ or $-I$.   
\end{description} 

The proof of Proposition \ref{prop: homlifting} which allows us to lift a positive homotopy of generic loops fails if a loop is non-generic.  To connect the ${H}^i(s_{i+1}, t)$ to
${H}^{i+1}(s_{i+1},t)$ via positive loops using the Proposition \ref{prop: Dusalifting} , we need to know that $h(s_{i+1},t)$ is a generic loop in ${\Co}$.
However, the argument can be patched rather easily for the particular homotopy
$h(s,t)$ constructed above.  It is enough to show how to produce ${H}$ locally around $s_{i+1}$ when $h(s_{i+1},t)$ has one
diversion into  $\pi({\Bb}_{\Uu , \Dd})$ or ${\pi}({\Bb}_{\Rr, \Dd})$
or $N_1^{-,-}$ as produced in Lemma 
\ref{lemma: outofOc}  and when
there are finitely many points at $I$ or $-I$ as in Lemma 
\ref{lemma: 4x4std}.  The final three lemmas  complete our discussion.

\begin{lemma}
If $h(s_{i+1},t)$ is non generic because it enters
${\pi}({\Bb}_{\Uu, \Dd})$ at time $t = t_0$ as in Lemma \ref{lemma: type(1)},
we can construct a local lifting of $h$.
\end{lemma}

\proof{}
 In the proof of Lemma \ref{lemma: type(1)} we actually constructed
a lift  of $h$ for $s,t$ in some interval $[s_{i+1} - \epsilon,
s_{i+1} + \epsilon] \times [t_0 - \delta, t_0 + \delta]$.  However,
the paths at $s = s_{i+1} \pm \epsilon$ are not generic, as they still go
through $\Bb_{\Uu, \Dd}$ at time $t = t_0$.  It is not hard to see
that one can stil patch these different local lifts by the argument of
Proposition \ref{prop: Dusalifting}. The important thing is that the
fibres of $\pi$ are always connected and there is only one non-generic
point on each path. 
 $\Box$

\begin{lemma}
If $h(s_{i+1},t)$ is non generic because it enters
${\pi}({\Bb}_{\Rr, \Dd})$ or $N_1^{-,-}$ as in Lemma \ref{lemma: outofOc},
we can construct a local lifting of $h$.
\end{lemma}

\proof{}
In the proof of this Lemma \ref{lemma: outofOc}, we
actually constructed a lift ${H}(s,t)$ of $h(s,t)$  for $s \in 
[s_{i+1} - \delta, s_{i+1} + \delta] $ for some $\delta > 0$ such that
$h(s_{i+1} - \delta,t)$  and $h(s_{i+1} + \delta, t)$  are generic loops
in ${\Co}$.    We can relabel
the $s_i$ appopriately and apply the remainder of the proof
of Proposition \ref{prop: homlifting} to lift the entire homotopy. $\Box$

\begin{lemma}
If $h(s_{i+1},t)$ is non generic because it 
passes through $I$ or $-I$ at times other than $0$ and $2 \pi$, we can
construct a local lifting of $h$. 
\end{lemma}
\proof{}
 By compactness,
there are finitely many such times, say $ \{t_j \} \mid _{1 \leq j \leq N} $.
Then, for each interval $[t_j, t_{j+1}]$, $h(s_{i+1},t)$ is is a positive
generic path in ${\Co}$ starting and ending at $I$ or $-I$.  Call this path
$h_j (s_{i+1},t)$.  By Lemma \ref{prop: Dusalifting}, the space of
positive lifts of $h_j (s_{i+1},t)$ is path connected.  Thus,  we can connect
${H}^i _j (s_{i+1}, t)$ to ${H}^{i+1} _j (s_{i+1}, t)$  for each
$1 \leq j \leq N$ independently, and arrive at a piecewise positive homotopy
in ${\Sp}(4)$ between ${H}^i (s_{i+1}, t)$ and ${H}^{i+1} (s_{i+1}, t)$.  Since piecewise positive paths can be approximated arbitrarily closely
by positive paths, we can find a positive homotopy in ${\Sp}(4)$ between
${H}^i (s_{i+1}, t)$ and ${H}^{i+1} (s_{i+1}, t)$.  As in the 
proof of Proposition \ref{prop: homlifting}, we patch together the ${H}^i (s_{i+1}, t)$ and ${H}^{i+1} (s_{i+1}, t)$ to obtain ${H}(s,t)$. $\Box$

\section{Technical Proofs}

This section contains the proofs of the technical lemmas needed in
Section 4.  We will restate them here for the convenience of the
reader.

\begin{makingDe}
 There exists a path $\De ^s : [0,1] \rightarrow  {\Bb}_{\Rr} \cup {\Bb}_{\Uu}^- \cup {\Nn}_1^{-,-}$ where $\De^0 \in {\Bb}_{\Rr}$, $\De^{\frac {1}{2}}=
\Nn_1^{-,-}$, and $\De^1 \in {\Bb}_{\Uu}^-$ satisfying two properties:
\begin{description}
\item[(i)] There exists a (positive) vector field along $\De ^s$ pointing into
${\Oo}_{\Cc}$  of the form $JP{\De}^s$ for a positive definite $P$.
\item[(ii)] There exists a negative vector field along $\De ^s$ pointing
into ${\Oo}_{\Rr}^+$ when $\De ^s \in {\Bb}_{\Rr}$ and pointing into
${\Oo}_{\Uu}$ elsewhere.  
\end{description}
\end{makingDe}
\proof{}

First, we need to construct the path ${\De}^s$.  The first part of ${\De}^s$ will travel
within the boundary components from ${\De}^0 \in {\Bb}_{\Rr} $ where
$\pi({\De}^0) = {\de}^0$ to 
${\De}^{\epsilon} \in {\Bb}_{\Rr, \Dd}$.  

Suppose that ${\de}^{\epsilon} = \pi(\Delta^\epsilon) \in {\pi}({\Bb}_{\Rr, \Dd}) $ has eigenvalues $ {\la}, {\la}, \frac{1}{\la}, \frac{1}{\la}$ and ${\de}^1 \in {\pi}({\Bb}_{\Uu}^-)$ has
eigenvalues $ a+bi, a+bi, a-bi, a-bi$ where $ a^2 + b^2 = 1$.  
Let ${\De}^s : [0,1] \rightarrow {\Sp} (4)$ be the path defined as 
$$ {\De}^s = \left \{ \begin{array} {lll} 
\left ( \begin{array} {cccc} 
\mu & 1 & 0 & 0 \\
0 & \frac{1}{\mu} & 0 & 0 \\
0 & 0 &  \mu & 1 \\
0 & 0 & 0 & \frac{1}{\mu} \end{array} \right ) & if & 
\epsilon < s < \frac{1}{2} \\
\left ( \begin{array} {cccc}
1 & 1 & 0 & 0 \\
0 & 1 & 0 & 0 \\
0 & 0 & 1 & 1 \\
0 & 0 & 0 & 1 
\end{array} \right ) & if & s = \frac{1}{2} \\
\left ( \begin{array} {rrrr}
x & x & \sqrt{1-x^2} & \sqrt{1-x^2} \\
0 & x & 0 & \sqrt{1-x^2} \\
- \sqrt{1-x^2} & - \sqrt{1-x^2} & x &  x \\
0 & - \sqrt{1-x^2} & 0 & x 
\end{array} \right ) & if & s > \frac{1}{2}
\end{array} \right .$$
where $x = 2-a+(2a-2)s$ and $\mu = \frac{ {\la}-1} {\epsilon - \frac{1}{2}} (s - \frac{1}{2})+1$.
Then ${\pi}({\De}^s) = {\de}^s$ lies in the appropriate regions.

We will now look for a positive continuous vector field along ${\De}^s$ which points into ${\Oo}_{\Cc}$ at every point
and project it to ${\Co}$ to get
the needed vector fields along ${\de}^s$.  The original path $A_t$ gives us one positive vector, say $v_0$, pointing into ${\Oo}_{\Cc}$ at ${\De}^0$.  We claim that $JP {\De}^s$ is a positive vector at ${\De}^s$ pointing into ${\Oo}_{\Cc}$ for all $\epsilon < s \leq 1$, where 
$$ P = \left ( \begin{array} {rrrr} 
10 & 0 & 0 & 1 \\
0 & 10 & 0 & 0 \\
0 & 0 & 10 & 0 \\
1 & 0 & 0 & 12 \\
\end{array} \right ) .$$ 
Since the positive cone is open and convex, join $v_0$ to $JP\Delta^\epsilon$
by a family of poisitive vectors pointing into ${\Oo}_{\Cc}$ along the
path $\Delta^s$ for $0<s<\epsilon$.    Then, we can continue the
vector field along ${\De}^s$  by letting the tangent vector at time $s$ equal
$JP {\De}^s$ for all $\epsilon \leq s \leq 1$.  This vector field is certainly continuous and positive, we need only prove the claim that it points
into  ${\Oo}_{\Cc}$ for all time.

When $s > \frac{1}{2}$, ${\De}^s \in {\Bb}_{\Uu}^-$, and thus any positive
vector points into ${\Oo}_{\Cc}$.  Also, by construction, our positive vector
field points into  ${\Oo}_{\Cc}$ for $s < {\eps}$.   Hence, we need only consider ${\eps} \leq s \leq \frac{1}{2}$.  We check the direction of these vectors by examining the behavior of the symmetric functions of the eigenvalues of paths in their directions.  For all matrices in  ${\Bb}_{\Rr , \Dd} \cup {\Bb}_{\Rr} \cup {\Nn}_1 ^{-,-}$, ${\si}_2 = \frac{{\si}_1 ^2}{4} + 2$ while, on the other hand, matrices in ${\Oo}_{\Cc}$
satisfy   ${\si}_2 > \frac{{\si}_1 ^2}{4} + 2$ and those in ${\Oo}_{\Rr}$ satisfy ${\si}_2 < \frac{{\si}_1 ^2}{4} + 2$.

We look at the derivatives 
$$\begin{array} {c}
\frac{d}{dr} |_{r=0} {\si}_1 (e^{JPr} {\De^s}) = {\si}_1 ' (s) \\
\frac{d}{dr} |_{r=0} {\si}_2 (e^{JPr} {\De^s}) = {\si}_2 ' (s)  \end{array}$$

Since ${\si}_2 = \frac{{\si}_1 ^2}{4} + 2$ for all points on ${\De}^s$ for
$s \leq \frac{1}{2}$, 
if ${\si}_2 ' (s) > (\frac{ {\si}_1 (s) ^2 }{4} +2)' $, then we know that
$JP{\De}^s$ points into ${\Oo}_{\Cc}$.  More generally, if $$
\frac{d^k}{dr^k} |_{r=0} ({\si}_2 (s)) = \frac{d^k}{dr^k} |_{r=0} (\frac{{\si}_1^2 (t)}{4} +2)$$
for all $k \leq n$, and $$
\frac{d^n}{dr^n} |_{r=0} ({\si}_2 (s)) > \frac{d^n}{dr^n} |_{r=0} (\frac{{\si}_1^2 (s)}{4} +2)$$
then $JP{\De}^s$ points into ${\Oo}_{\Cc}$.

Let us consider specifically the point ${\De}^\frac{1}{2} = {\Nn}_1^{-,-}$.
If we examine the symmetric functions of $e^{JQr} {\Nn}_1^{-,-}$ for general
symmetric $$
Q = \left ( \begin{array} {cccc}
q_1 & q_2 & q_3 & q_4 \\
q_2 & q_5 & q_6 & q_7 \\
q_3 & q_6 & q_8 & q_9 \\
q_4 & q_7 & q_9 & q_{10}  \end{array} \right ) $$
we find that ${\si}_2 ' = (\frac{ {\si}_1 ^2}{4} +2)'$ for all $Q$.  Going to the second derivative, ${\si}_2 '' < (\frac{{\si}_1 ^2}{4} +2)''$, except if $q_3 = 0$ and $q_1 = q_8$, in which case ${\si}_2 '' = (\frac{{\si}_1 ^2}{4} +2)''$. Imposing these two restrictions on Q and looking at the third derivatives, we find  ${\si}_2 ''' > (\frac{{\si}_1 ^2}{4} +2)'''$ if $q_1 > 0$ and  $q_4 \neq q_6$.  Hence, $e^{JQr} {\Nn}_1^{-,-}$ is a positive path
pointing into ${\Oo}_{\Cc}$, if $Q$ is a positive definite matrix 
satisfying $q_1 > 0$, $q_3 = 0$, $q_4 \neq q_6$, and $q_1 = q_8$.  Indeed, the
aforementioned matrix $P$ satisfies these conditions, and we can check that$$
\begin{array} {r}
{\si}_2 ' = (\frac{ {\si}_1 ^2}{4} +2)' = 20 \\
{\si}_2 '' = (\frac{ {\si}_1 ^2}{4} +2)'' = -680 \\
{\si}_2 ''' = -17560 \\
(\frac{ {\si}_1 ^2}{4} +2)''' = -17600 \end{array} $$
for the path $e^{JPr} {\Nn}_1^{-,-}$ , and hence this path does travel into
${\Oo}_{\Cc}$. 

Additionally, consider the path $e^{JPr}{\Nn}_{1+y}^{-,-}$ where $$
{\Nn}_{1+y} ^ {-,-} = \left ( \begin{array} {rrrr}
1+y & 1 & 0 & 0 \\
0 & \frac{1}{1+y} & 0 & 0 \\
0 & 0 &  1+y & 1 \\
0 & 0 & 0 & \frac{1}{1+y} \end{array} \right ) .$$
This path satisfies $$
\begin{array} {r}
{\si}_2 ' = (\frac{ {\si}_1 ^2}{4} +2)' \\
{\si}_2 '' > (\frac{ {\si}_1 ^2}{4} +2)'' \end{array} $$
for all $y > 0$.

The matrices in ${\De}^s$ for $ {\eps}< s< \frac{1}{2}$ are all of the form ${\Nn}_{1+y} ^ {-,-}$ for some  $y > 0$.  Therefore, the positive vector field which we have constructed on this portion of the path, $JP{\Nn}_{1+y}^{-,-}$ points into ${\Oo}_{\Cc}$ and the proof of the claim is
completed.

Finally, we need to construct a negative (so the reverse flow would be positive) vector field along ${\De}^s$ which points into
${\Oo}_{\Rr}$ in the direction of decreasing trace for $s< \frac{1}{2}$ and into ${\Oo}_{\Uu}$ for $s \geq \frac{1}{2}$.  For  $s \geq \frac{1}{2}$, ${\De}^s \in {\Bb}_{\Uu}^- \cup {\Nn}_1^{-,-}$, and 
all negative vectors based at  ${\De}^s$ will point into ${\Oo}_{\Uu}$.  
Therefore, if we find any negative continuous vector field along ${\De}^s$ for
$s< \frac{1}{2}$, any negative continuous extension of it will provide us with
vectors pointing into ${\Oo}_{\Uu}$ for the duration of ${\De}^s$.  We can pick
such an extension to match the tangent vector of ${\ga}$ at the point ${\De}^1$.

For $\epsilon \leq s \leq \frac{1}{2}$, on ${\De}^s$, we have block matrices of the form $$
\left ( \begin{array} {rrrr}
{\mu} & 1 & 0 & 0 \\
0 & \frac{1}{ {\mu}} & 0 & 0\\
0 & 0 & {\mu} & 1 \\
0 & 0 & 0 & \frac{1}{ {\mu}}   \end{array} \right ) . $$
It would be sufficient, then, to find a negative definite $2 \times 2$ matrix
$Q_2$  such that $$
JQ_2 \left ( \begin{array} {rr} {\mu} & 1 \\ 0 & \frac{1}{{\mu}} \end{array} \right ) $$
points into ${\Oo}_{\Rr}$ in the direction of decreasing trace for all ${\mu}$.  Then, set
$Q_4$ equal to the $4 \times 4$ block matrix with $Q_2$ in the upper left
and lower right blocks, and vector field $JQ_4 {\De}^s$ is a negative, continuous vector field pointing into ${\Oo}_{\Rr}$ in the direction of decreasing trace for $\epsilon \leq s \leq \frac{1}{2}$.  For $s < \epsilon$, we can continuously perturb $Q_4$ so that $JQ_4{\De}^s$ is a negative vector field pointng into ${\Oo}_{\Rr}$ in the direction of decreasing trace which matches the 
given tangent vector to ${A_t}$ at ${\De}^0$.  However, matrices $Q_2$ are plentiful;
one can be chosen which can be slightly perturbed along ${\De}^s$ to match the tangent vector
to ${A_t}$ at ${\De}^0$. $\Box$

\begin{makingdifforder}
Let $k,\ell \geq 1$, $n \geq 3$, $n>k$, $n>\ell$ where $k,\ell,n$ are all integers.
Then, $$
\left ( \begin{array} {rr}
e^{Jkt} & 0 \\
0 & e^{J(n-k)t} \end{array} \right ) \sim_+
\left ( \begin{array} {rr}
e^{J\ell t} & 0 \\
0 & e^{J(n- \ell )t} \end{array} \right ) $$
where $J \in {\Sp}(2)$.
\end{makingdifforder}

\proof{ }
 
The positive homotopy between the two paths is $$
H({\theta}, t) = \left ( \begin{array} {rr}
e^{Jt} & 0 \\
0 & e^{J(1+n-k- \ell )t} \end{array} \right ) P_{\theta} 
\left ( \begin{array} {rr}
e^{J(k-1)t} & 0 \\
0 & e^{J(\ell -1)t} \end{array} \right ) P_{\theta}^{-1} $$
for $t \in [0, 2{\pi}]$ and ${\theta} \in [0, \frac{{\pi}}{2}]$ where $$
P_{\theta} = \left ( \begin{array} {rr}
\cos ({\theta})I & - \sin ({\theta})I \\
\sin({\theta})I & \cos({\theta}) I \end{array} \right ) \in {\Sp}(4). $$
Here, $I$ represents the $2 \times 2$ identity matrix.
$H({\theta}, t)$ is certainly a homotopy, as it is the product of symplectic
matrices for all time and hence always contained in ${\Sp}(4)$, and $$
H(0,t) = \left ( \begin{array} {rr}
e^{Jkt} & 0 \\
0 & e^{J(n-k)t} \end{array} \right ) $$
$$ H(\frac{{\pi}}{2}, t) = \left ( \begin{array} {rr}
e^{J \ell t} & 0 \\
0 & e^{J(n-\ell )t} \end{array} \right ) . $$
We must check that this is a positive homotopy, i.e. $H({\theta}, t)$ is a positive path for any fixed ${\theta} \in [0, \frac{{\pi}}{2} ]$.  
Let $R$ be the $4 \times 4$ matrix such that  $$
\frac{d}{dt} |_{t=t_0} H({\theta}, t) = JRH({\theta}, t_0). $$
Certainly, $R$ depends on both ${\theta}$ and $t_0$.   $H({\theta}, t)$
is positive if and only if $R$ is a positive definite matrix for all ${\theta}$
and for all $t_0$.  $R$ must be
symmetric since $JRH({\theta}, t_0)$ is in the tangent space of $ {\Sp}(4)$
at the point $H({\theta}, t_0)$, thus it will be sufficient to prove that
the eigenvalues of $R$ are positive real.

Without loss of generality, assume that $k > \ell $ and $k, \ell \leq \frac{n}{2}$.
The second assumption is justified because , $$
\left ( \begin{array} {rr}
e^{Jkt} & 0 \\
0 & e^{J \ell t} \end{array} \right ) \sim_+
\left ( \begin{array} {rr}
e^{J \ell t} & 0 \\
0 & e^{Jkt} \end{array} \right ) $$
under the positive homotopy $$
G({\theta},t) = P_{\theta} \left ( \begin{array} {rr}
e^{Jkt} & 0 \\
0 & e^{J \ell t} \end{array} \right ) P_{\theta}^{-1} $$
for ${\theta} \in [0, \frac{{\pi}}{2}]$.  
$G({\theta},t)$ is positive for any fixed ${\theta}$ since it is the
conjugate of a positive path, and $$
G(0,t) = \left ( \begin{array} {rr}
e^{Jkt} & 0 \\
0 & e^{J \ell t} \end{array} \right )$$
$$ G(\frac{{\pi}}{2},t) = \left ( \begin{array} {rr}
e^{J \ell t} & 0 \\
0 & e^{Jkt} \end{array} \right ) . $$

We now compute $R$ to determine its eigenvalues.
Let $J$ denote both the standard $2 \times 2$ and $4 \times 4$ matrix, its
dimension will be clear by context.  Let $r = 1+n-k- \ell $ to make computations
easier. 
$$\begin{array} {rl} \frac{d}{dt} H({\theta},t) =  &
\left ( \begin{array} {rr} Je^{Jt} & 0 \\
                            0 & rJe^{Jrt} \end{array} \right )
P_{\theta} \left ( \begin{array} {rr} e^{J(k-1)t} & 0 \\
                         0 & e^{J(\ell -1)t} \end{array} \right ) P_{\theta}^{-1} + \\ & \\
 & \left ( \begin{array} {rr} e^{Jt} & 0 \\
                            0 & e^{Jrt} \end{array} \right ) P_{\theta}
\left ( \begin{array} {rr} (k-1)Je^{J(k-1)t} & 0 \\
                            0 & (\ell -1)Je^{J(\ell -1)t} \end{array} \right ) P_{\theta} ^{-1} \\ & \\ & \\
= & J \left ( \left ( \begin{array} {rr} I & 0 \\ 0 & rI \end{array} \right ) 
 +  J^{-1} \left ( \begin{array} {rr} e^{Jt} & 0 \\
                            0 & e^{Jrt} \end{array} \right ) P_{\theta} \times \right. \\ & \\

& \left. \left ( \begin{array} {rr} (k-1)Je^{J(k-1)t} & 0 \\
                            0 & (\ell -1)Je^{J(\ell -1)t} \end{array} \right ) P_{\theta}^{-1} H({\theta},t)^{-1} \right)  H({\theta},t) \\
\end{array} . $$

Multiplying the terms in the parentheses gives $$
R
=  \left ( \begin{array} {rr} 
(1+(k-1) \cos ^2{\theta} + (l-1) \sin ^2{\theta})I & 
\cos {\theta} \sin {\theta} (k-\ell) e^{Jt(1-r)} \\
\cos {\theta} \sin {\theta} (k-\ell) e^{Jt(r-1)} &
(r+(k-1) \sin ^2 {\theta} + (\ell-1) \cos ^2 {\theta})I \end{array} \right ) $$
$R$ has two eigenvalues of multiplicity two which happen to be independent of $t$:
$$ {\la}_1 = \frac{1}{2} (n + \sqrt{ (k-\ell)^2 + 2 \cos (2{\theta})(k-\ell)(1-r) + (1-r)^2}) $$
$$ {\la}_2 = \frac{1}{2} (n - \sqrt{ (k-\ell)^2 + 2 \cos (2{\theta})(k-\ell)(1-r) + (1-r)^2}) . $$

Certainly, since $n$ is positive, ${\la}_1$ is positive for all ${\theta}$.  To check that 
${\la}_2$ is positive, we must show $$
\begin{array} {rrr} 
\sqrt{ (k-\ell)^2 + 2 \cos (2{\theta})(k-\ell)(1-r) + (1-r)^2} & < & n \end{array}$$
Recall the previously justified assumptions that $k>\ell$ and $k,\ell \leq \frac{n}{2}$.  If $k=\ell=\frac{n}{2}$, then $r=1$ and the left hand side of the inequality
is $0$ which is certainly less than $n$. If, on the other hand, either $k$ or
$\ell$ is less than $\frac{n}{2}$, then $(1-r)$ is negative while $(k-\ell)$ is positive.  Hence, $$
\begin{array} {rcl}
 \sqrt{ (k-\ell)^2 + 2 \cos (2{\theta})(k-\ell)(1-r) + (1-r)^2} & \leq &
\sqrt{(k-\ell)^2 -2(k-\ell)(1-r) + (1-r)^2} \\
& = & \sqrt{ ((k-\ell) - (1-r))^2} \\
& = & k-\ell-1+r  \\
& = & n-2l \\
& < & n \end{array} $$
and thus ${\la}_2$ is positive for all ${\theta}$.  Hence, $R$ is a positive
definite matrix, and $H({\theta},t)$ is a positive homotopy. $\Box$

\end{document}